\documentclass[10pt]{amsart}
\usepackage{color}
\usepackage{verbatim}

\usepackage[square,compress,comma, numbers,sort]{natbib}
\usepackage[colorlinks=true, citecolor=green, linkcolor=blue]{hyperref}

\allowdisplaybreaks[4]

\def\Var{\text{Var}}
\def\H{\mathcal{H}}

\def\ooo{(1+o(1))}

\def\LT{\left}
\def\RT{\right}

\definecolor{c20}{rgb}{0.,0.7,0.}
\definecolor{c30}{rgb}{0.,0.,1.}
\definecolor{c40}{rgb}{1,0.1,0.7}
\definecolor{c50}{rgb}{1,0,0}
\definecolor{c60}{rgb}{1,0.9,0.1}

\def\b{\vk{b}}
\def\x{\vk{x}}
\def\tilb{\widetilde{\b}}
\def\tilbu{\widetilde{\b_u}}

\def\SI{\Sigma}

\newcommand{\abs}[1]{\left\lvert #1 \right\rvert}

\newcommand{\E}[1]{\mathbb{E}\left\{ #1\right\}}

\newcommand{\pk}[1]{\mathbb{P} \left\{ #1 \right \} }

\newcommand{\R}{\mathbb{R}}

\newcommand{\N}{\mathbb{N}}
\newcommand{\inr}{\in \R}

\newcommand{\ldot}{,\ldots,}

\newcommand{\limit}[1]{\lim_{#1 \to   \infty}}

\newcommand{\BQN}{\begin{eqnarray}}
\newcommand{\EQN}{\end{eqnarray}}
\newcommand{\BQNY}{\begin{eqnarray*}}
\newcommand{\EQNY}{\end{eqnarray*}}

\newcommand{\BS}{\begin{sat}}
\newcommand{\ES}{\end{sat}}
\newcommand{\BT}{\begin{theo}}
\newcommand{\ET}{\end{theo}}
\newcommand{\BK}{\begin{korr}}
\newcommand{\EK}{\end{korr}}
\newcommand{\EQD}{\stackrel{d}{=}}

\newcommand{\BD}{\begin{de}}
\newcommand{\ED}{\end{de}}

\newcommand{\BIT}{\begin{itemize}}
\newcommand{\EIT}{\end{itemize}}
\newcommand{\BDI}{\begin{description}}
\newcommand{\EDI}{\end{description}}

\newcommand{\BRM}{\begin{remarks}}
\newcommand{\ERM}{\end{remarks}}

\newcommand{\BEL}{\begin{lem}}
\newcommand{\EEL}{\end{lem}}

\newtheorem{theo}{Theorem}[section]
\newtheorem{sat}[theo]{Proposition}
\newtheorem{de}[theo]{Definition}
\newtheorem{lem}[theo]{Lemma}

\newtheorem{korr}[theo]{Corollary}
\newtheorem{remark}[theo]{Remark}
\newtheorem{remarks}[theo]{Remarks}

\newcommand{\nelem}[1]{{Lemma \ref{#1}}}

\newcommand{\netheo}[1]{{Theorem \ref{#1}}}
\newcommand{\nekorr}[1]{{Corollary \ref{#1}}}

\newcommand{\prooftheo}[1]{ \textsc{\bf Proof of Theorem} \ref{#1}:}

\newcommand{\prooflem}[1]{\textsc{\bf Proof of Lemma} \ref{#1}:}

\newcommand{\COM}[1]{}

\newcommand{\QED}{\hfill $\Box$}

\newcommand{\kb}[1]{\boldsymbol{#1}}
\newcommand{\vk}[1]{\kb{#1}}

\def\IF{\infty}

\def\LT{\left}
\def\RT{\right}

\def\ooo{(1+o(1))}

\def\vn{\varepsilon}

\def\Del{\triangle}

\def\uY { {\overline{\vk{u}}}}

\def\tDu{\Delta_u(\tau)}

\def\HAS{\mathcal{H}_I}

\def\X{\vk{X}}

\def\calU{\mathcal{U}}

\def\KJ{\mathcal{K}_J}
\def\tDu{\Delta_u(\lambda,\tau)}

\def\xmu{ (\X (t)-\vk{\mu}t) }
 \def\squ {\sqrt{u}}

\def\K1#1{\textcolor{black}{#1}}
\def\cc#1{\textcolor{black}{#1}}
\def\EH#1{\textcolor{black}{#1}}

\def\ccj#1{\textcolor{black}{#1}}

\def\EHb#1{\textcolor{black}{#1}}

\def\K1#1{\textcolor{black}{#1}}
\def\EHc#1{\textcolor{black}{#1}}
\def\Lc#1{\textcolor{black}{#1}}

\def\k1#1{\textcolor{black}{#1}}

\def\Jc#1{\textcolor{black}{#1}}
\def\e1#1{\textcolor{black}{#1}}
\def\kk#1{\textcolor{black}{#1}}

\def\LL#1{\textcolor{black}{#1}}

\def\wHAS{\widetilde{\HAS}}

\def\ci{C_I}
\def\gt{\widehat g}
\def\ggt{ \widetilde{ g}}

\begin{document}

\title{Extremal behaviour of hitting a cone by correlated Brownian motion with drift}

\author{Krzysztof D\c{e}bicki}
\address{Krzysztof D\c{e}bicki, Mathematical Institute, University of Wroc\l aw, pl. Grunwaldzki 2/4, 50-384 Wroc\l aw, Poland}
\email{Krzysztof.Debicki@math.uni.wroc.pl}

\author{Enkelejd  Hashorva}
\address{Enkelejd Hashorva, Department of Actuarial Science, 
University of Lausanne,\\
UNIL-Dorigny, 1015 Lausanne, Switzerland
}
\email{Enkelejd.Hashorva@unil.ch}

\author{Lanpeng Ji}
\address{Lanpeng Ji, Department of Actuarial Science, 
University of Lausanne\\
UNIL-Dorigny, 1015 Lausanne, Switzerland
}
\email{jilanpeng@126.com}

\author{Tomasz Rolski}
\address{Tomasz Rolski, Mathematical Institute, University of Wroc\l aw, pl. Grunwaldzki 2/4, 50-384 Wroc\l aw, Poland}
\email{Tomasz.Rolski@math.uni.wroc.pl}

\bigskip

\date{\today}
 \maketitle

 {\bf Abstract:}
This paper derives an exact asymptotic expression for
\[
\mathbb{P}_{\mathbf{x}_u}\{\exists_{t\ge0}  \mathbf{X}(t)-  \boldsymbol{\mu}t\in \mathcal{U} \}, \ \ {\rm as}\ \ u\to\infty,
\]
where
$\mathbf{X}(t)=(X_1(t),\ldots,X_d(t))^\top,t\ge0$ is a correlated $d$-dimensional Brownian motion
starting at the point
$\mathbf{x}_u=-\boldsymbol{\alpha}u$ with $\boldsymbol{\alpha}\in \mathbb{R}^d$,
$\boldsymbol{\mu} \in \mathbb{R}^d$ and $\mathcal{U}=\prod_{i=1}^d [0,\infty)$.
The derived asymptotics depends on the
solution of an underlying multidimensional quadratic optimization problem with constraints, which
leads in some cases to dimension-reduction of the considered problem. Complementary, we study
asymptotic distribution of the conditional
first passage time to $\mathcal{U}$,
which depends on the dimension-reduction phenomena.

 {\bf Key Words:} multidimensional Brownian motion; extremes; exact asymptotics;  \EHb{first passage time}; large deviations; quadratic programming problem;  multidimensional Pickands constants.

 {\bf AMS Classification:} Primary 60G15; secondary 60G70

 \section{Introduction}
Consider $\vk{X}(t)-\vk{\mu}t , t\ge0$, a correlated $d$-dimensional Brownian motion with drift, where $\vk{X}(t)=A\vk{B}(t)$,
$A \in \R^{d\times d}$ is a non-singular matrix,
$\vk{B}(t)=(B_1(t),\ldots,B_d(t))^\top,t\ge0$ is a standard $d$-dimensional Brownian motion with independent coordinates
and $\vk{\mu}=(\mu_1 \ldot \mu_d)^\top\cc{\in \R^d}$.

\k1{The probability
\begin{eqnarray}\label{e.0}
\mathbb{P}_{\vk{x}}\{\exists_{t\ge0}  \vk{X}(t)-  \vk{\mu}t\in \calU \}
\end{eqnarray}
that
starting at the point
$\vk{x}\in \R^d$, the process $\vk{X}(t)-\vk{\mu}t$
enters the set $\calU\subset\R^d$ in a finite time,
is of interest both for theory-oriented studies and for
applied-mathematics problems as, e.g.,
heat and mass diffusion, photon absorption or chemotaxis.
Due to the complexity of
(\ref{e.0}), still only some fragmentary results
focusing on the special case of mutually independent coordinates
(i.e., for $A$ being the identity matrix) or on particular structures of $\calU$
are available.
We refer to, e.g., \cite{KnK00}
for the asymptotic analysis, as $r:=\|\vk{x}\|\to\infty$, of
(\ref{e.0}) for $A$ the identity matrix,
some compact $\calU$, and appropriately chosen drifts, see also \cite{LLSB84,BSU}.
Somehow related problem for the exit time from a cone for
a (noncorrelated) multidimensional Brownian motion with drift
was considered in \cite{GaR14} and references therein;
see also \cite{PuR08}
for the case of $\calU$ being a Weyl chamber.\\
This contribution is concerned with
investigation of
(\ref{e.0}) for the model allowing
correlation between the Brownian components.
More precisely, we
investigate the \Jc{asymptotics of} probability that in infinite-time horizon,
the process
$\vk{X}(t)-\vk{\mu}t, t\ge0$,
starting at point
$\vk{x}_u:=(-\alpha_1u,...,-\alpha_d u)^\top$
\Lc{with $\alpha_i\in\R, 1\le i\le d, u>0$},
enters  the cone $\calU=\prod_{i=1}^d [0,\infty)$, that is}
\BQN \label{eq:Pu}
P(u):=
\mathbb{P}_{\vk{x}_u}\{\exists_{t\ge0}  \vk{X}(t)-  \vk{\mu}t\in \calU \}, \ \ \ \Jc{u\to\IF}.
\EQN
\k1{
Our results allow  for considering other sets in (\ref{eq:Pu}), as e.g.,
polyhedral cones
$\{\vk{x}\in \R^d: M\vk{x}\ge\vk{0}\}$, where \Jc{$M$} is a  \Jc{$d\times d$ non-singular} matrix.
Indeed, by a linear transformation of $M$,
we can reduce the problem of hitting the polyhedral cone \e1{to (\ref{eq:Pu}), namely}
\[
\mathbb{P}_{\vk{x}_u}\left\{\exists_{t\ge0}  \vk{X}(t)-  \vk{\mu}t \in \{\vk{x}\in \R^d: M\vk{x}\ge\vk{0}\} \right\}
=
\mathbb{P}_{\vk{x}^{,}_u}\left\{\exists_{t\ge0}  M\vk{X}(t)-  M\vk{\mu}t\in \prod_{i=1}^d [0,\infty) \right\},
\]
with $\vk{x}^{,}_u=M\vk{x}_u$.
}

\K1{Since we are} interested in the case that
\k1{$\limit{u}P(u)= 0$} we shall assume that there exists some $1\le i\le d$ such that
\BQN \label{alphamu}
	\alpha_i>0,    \quad  \mu_i>0.
\EQN
\k1{Using that
\[
P(u)=\mathbb{P}_{\mathbf{0}}
\LT\{\exists_ {t\ge 0}  \bigcap_{i=1}^{ d} \{ X_i(t)-  \mu_it>  \alpha_i u\}\RT\}
\]
this paper contributes also to extreme value problems of vector-valued stochastic processes.}

\K1{Complementary, we investigate distributional properties of the passage time
of $\vk X(t)-\vk \mu t$ to
\k1{$\calU$, for $\|\vk{x}_u\|\to\infty$ as $u\to\infty$,}
given that the multivariate process has ever entered the upper quadrant.}
Specifically, for
	\BQN \label{eq:RT}
	\tau_u=\inf\{t\ge 0: \vk X(t)-\vk \mu t>\vk\alpha u\}
	\EQN
\k1{($\mathbf{X}(0)=\vk{0}$)}
we are interested in the approximate distribution of $\tau_u \lvert  \tau_u<\IF$ as $u\to \IF$.

In the 1-dimensional \EHb{setup} it is well-known that for $\alpha, \mu$ positive
$$P(u)=\pk{\sup_{t\ge0} (B_1(t)-  \mu t)>  \alpha u}=  e^{-2 \alpha \mu u  }, 
$$
\k1{where from this point on we write
$\mathbb{P}:=\mathbb{P}_{\mathbf{0}}$.}
Further, in view of \cite{MR2462285} we have that
$$ \limit{u} \pk{ \cc{\alpha^{-1/2}}\mu^{3/2}  (\tau_u - \cc{\alpha}u /\mu)/\sqrt{u} \le s \Bigl \lvert \tau_u < \IF} =\Phi(s), \quad s \inr,   $$ 	
with $\Phi$ the distribution \Jc{function} of an $\ccj{\mathcal N}(0,1)$ random variable.
\EHb{Normal or exponential
approximations for 1-dimensional Gaussian \K1{counterparts of the considered \Lc{model in this contribution} are discussed} in \cite{MR2462285,DHJParisian,DHJ13a}.}

In the case  $d\ge 2$,  both the  approximation of $P(u)$ and the \Jc{approximate} distribution of $\tau_u \lvert \tau_u < \IF$
depend on \EHb{the} solution of a \EHb{related} quadratic optimization problem.
\K1{In particular,}
in the light of \cite{RolskiSPA}[Theorem 1], the logarithmic
asymptotics of (\ref{eq:Pu}) can be derived and takes the following form (hereafter $\sim$ means asymptotic equivalence as $u\to \IF$)
\BQN \label{eq1D}
-\ln P(u) &\sim & \frac{\gt}{2} u, \quad \gt=\inf_{t\ge 0} g(t),
 \EQN
 with
\BQN \label{fgt}
g(t)&=&\frac{1}{t} \inf_{\vk{v} \ge \vk{\alpha}+\vk{\mu} t}  \vk{v}^\top \Sigma^{-1}  \vk{v}, \quad
\Sigma=AA^\top.
\EQN
Clearly, \eqref{eq1D} is of no use for the approximation of the conditional passage time $\tau_u \lvert \tau_u < \IF$ as $u\to \IF$.

Our main result presented in Theorem \ref{Thm1} shows that
\BQN  \label{e.main2}
P(u)  \sim  \ci\HAS u^{\frac{1- m}{2}}e^{- \frac{\gt}{2} u},
\EQN
where $\ci>0$, $m\in \N$  are known constants and $\HAS$ is a multidimensional counterpart of the
celebrated Pickands constant that appears in the
extreme value theory of Gaussian random fields; see e.g.,  \cite{PickandsA, PicandsB, Pit96, DRolski, Harper2, DM, Dancheng, SBK}. In \netheo{KorrRT}
we derive approximation of the conditional passage time.

One of \K1{the findings} of this paper is that the set of \ccj{indexes} $\{1,\ldots,d\}$ of \ccj{the} \EHb{vector-}process $\vk{X}$ can be partitioned into
 \EHb{three subsets $I,J,K$. The index set}  $I$ determines \ccj{$m,$} $\gt$ and  $\HAS$  in the asymptotics \eqref{e.main2},
 \EHb{whereas both  $I$ and $K$ determine the constant $\ci$}. \EHb{\K1{Moreover,} the set $J$, whenever non-empty, contains indices that do
 not play any role in our asymptotic consideration.}
 \K1{Interestingly, the limit distribution of the conditional passage time \k1{derived in \netheo{KorrRT}} is Gaussian only if $K=\emptyset$.}

\EHb{Our investigation shows that} for $d\ge 2$, the  problem (\ref{eq:Pu}) is surprisingly hard even for the seemingly simple case
of independent components, \K1{that is} with $A$ being the identity matrix.
Besides, solving this particular case does not reveal the essential ingredients that determine the asymptotics of $P(u)$ in the general case
where $\Jc{A}$ is not the identity matrix.

\K1{The strategy of the proof of the main result, given in Theorem \ref{Thm1}, although in its roots based on the
{\it double sum technique} developed in 1-dimensional setting for extremes of Gaussian processes and fields (see, e.g. \cite{PickandsA, PicandsB,Pit96}),
needed new ideas that in several key steps of the argumentation significantly differ from methods used in 1-dimensional case.
In particular, one of difficulties is the lack of Slepian-type inequalities that could be applied in our vector-valued setting.
Notice also that the standard techniques utilized for proving the negligibility of the double-sum, as  e.g., in  \cite{Pit96}, do not
work in the general $d$-dimensional vector-valued case.
Other difficulty lies in analysis of
the multidimensional Pickands constants $\HAS$.
Establishing its finiteness and positivity requires significant efforts.}
\k1{The developed in this paper approach opens some possibilities
for its application to
asymptotic analysis of some related functionals of vector-valued
Gaussian processes.}
\COM{This and some other  observation can be seen when analyzing in details (see Section \ref{ss.two-dim} with proofs in   \ref{ss.anal.2-dim}) a 2-dimensional setup:
 if $\alpha_1> \alpha_2 > 0$ and
$
\Sigma=\left(
\begin{array}{cc}
1 & \rho\\
\rho & 1 \\
\end{array}
\right)
.$}

\COM{
\begin{itemize}
\item[i)] $\HAS$  is a complicated constant which depends on some unique index set $I\subseteq \{ 1 \ldot d\}$. Establishing its finitness and positivity requires significant efforts and new ideas;
\item[ii)]  It is possible to have similar asymptotics of $P(u)$ for $d\ge 2$ as in the 1-dimensional case. For instance for the 2-dimensional setup, if $\alpha_1> \alpha_2 > 0$ and
$
\Sigma=\left(
\begin{array}{cc}
1 & \rho\\
\rho & 1 \\
\end{array}
\right)
,$ then for \cc{$\rho\in ((\alpha_1+ \alpha_2)/(2 \alpha_1),1)$} we have (see below \nekorr{cor:TD})
\BQN \label{bauarb}
 \pk{\exists_{t\ge 0} \cap_{i=1}^2 \{(X_i(t) - t) > \alpha_i u\}} \sim
\ci \pk{\exists_{t\ge 0} (X_1(t) - t) > \alpha_1 u},\ \ u\to\IF,
 \EQN
with  $\ci=1$,   and thus only the first component of $\vk X(t), t\ge 0$ is controlling the asymptotics of $P(u)$. This case will be referred to  as the {\it loss of dimensions phenomena}. There are several technical issues related to the loss of dimensions as it will be explained in our proofs below;
\item [iii)]  There are other cases of  {\it loss of dimensions phenomena}, where some components \cc{other than those with indexes in $I$} still play a role in the asymptotics of $P(u)$, but only up to some constants. For instance, referring again to the  2-dimensional case presented in \nekorr{cor:TD} we have for $\rho=(\alpha_1+ \alpha_2)/(2 \alpha_1)$ that
  \eqref{bauarb} holds, with \cc{$\ci$ taking the information of the second component and given by}
  $$ \ci=   \frac{1 }{\sqrt{ 2\pi \alpha_1}}\int_{\R}e^{-\frac{1}{2\alpha_1}x^2}\Psi\LT(\frac{1-\rho}{\sqrt{\alpha_1}}x\RT)\, dx,
  $$
where
  $\Psi(x)=1- \Phi(x)$ denotes  the survival function of an $\EH{\mathcal{N}}(0,1)$ random variable;
\item[iv)] The standard techniques utilized for proving the negligibility of the double-sum, see e.g., \cite{Pit96} do not
  work in the general $d$-dimensional vector-valued case of this paper. Here we make strong use of the independence of increments and the self-similarity of Brownian motion.
  \end{itemize}
}

\COM{
One of discoveries of this paper is that the index set of components $\{1,\ldots,d\}$ of the  process $\vk{X}(t)-\vk \mu t, t\ge0$ can be
partitioned into three subsets $I, J$ and $K$. \EH{The index set $I$ of the}   {\it essential} components, determines \EH{ both the constant}   $\inf_{t\ge 0} g(t)$ \EH{in the exponent of the asymptotics} and $\gamma$. \EH{The components \cc{with index in $ J$} do not play any role in the asymptotics and are therefore referred to as {\it inessential}, whereas the components with index in $ K$ determine only the
 prefactor so we referr to $K$ as the index set of \cc{\it weakly essential index set}. As we shall show below}
$\inf_{t\ge 0} g(t)$ can be expressed
explicitly in a form given in \eqref{eq:intr1}.
}

In this contribution we present a full general picture and a complete solution of the problem at hand by 
\Jc{developing new techniques building up on} asymptotic theory,
convex optimization and probability theory.
\EHb{Additionally,  we analyze in details \K1{some} special cases including  the case of independent components,
the homogeneous case when $\alpha_j=\alpha$ and $\mu_j=\mu$ for all $j$ and the case with negatively associated components.
Moreover, we discuss several interesting special cases when $d=2$.}

We organise the paper as follows. The next section fixes the notation and presents some preliminary \EHb{findings}.
The main results with examples are presented in Section 3, with detailed proof relegated to Section \ref{proofmain}.
Detailed analysis of the \K1{related} optimization problem and \EHb{some technical} proofs are displayed in Appendix.

\COM{In the setting of general Gaussian processes,
multidimensional extensions of the asymptotic properties of
supremas were proposed in
\cite{Piterbarg05} for $d=2$ and centered stationary components and \cite{RolskiSPA}
in a general setting of \cc{non-centered} Gaussian processes, for which it was possible to derive logarithmic asymptotics.
Recently exact asymptotics of
$\pk{\exists_ {t\in [0,T]} \min_{i=1 \ldot d}Y_i(t)>u}$
for $Y_i(t)$ being centered and mutually independent Gaussian processes that satisfy some regularity properties
for covariance functions were derived in \cite{DHJT15}, \cite{DELK}.
In this paper we allow \EH{for} dependence between coordinates, which \EH{renders} the
techniques used in \cite{DHJT15} or \cite{DELK} \EH{inapplicable here}.
The general strategy of the proof of (\ref{eq:PuAsym})
extends the technique of {\it double sum}, introduced by Pickands \cite{PickandsA} and developed by Piterbarg \cite{Pit96},
originally dedicated to one-dimensional case.
To get a multidimensional counterpart of this method, one has to cope with
lack of Slepian-type inequalities \EHb{(see \cite{GennaBorell,LedouxA,LedouxB} for the Gaussian case and \cite{GennaSlepian} for
for stable processes)} and
derivation of the multidimensional counterpart of Pickands lemma.
An additional difficulty is related with the analysis of the
properties of $\inf_{t\ge 0} g(t)$, which 
plays \cc{a} key role for the form of the asymptotics.
In Section \ref{sectggg} we investigate this problem \EH{utilising some results from}
  \cite{HA2005}.
}

\section{Preliminaries}
All vectors here are $d$-dimensional column vectors written in bold letters with  \K1{$d\ge2$}. For instance
$\vk{\alpha}=(\alpha_1 \ldot \alpha_d)^\top$, with $^\top$ the transpose sign.
Operations with vectors are meant component-wise, so 
$\abs{\vk{x}}=(\abs{x_1}   \ldot \abs{x_d })^\top$ and $\lambda \x=\x \lambda = ( \lambda x_1 \ldot \lambda x_d)^\top$ for any $\lambda\inr, \x\inr^d$. 
We denote
$$\vk{0}=(0\ldot 0)^\top \inr^d, \quad \vk{1}=(1\ldot 1)^\top \inr^d. $$
For any \EHb{non-empty subset} $\mathcal T\subset \R$, denote the inner set of $\mathcal T$ by $\mathcal T^o$ and its closure set by $\overline{\mathcal T}$.
If $I\subset \{1,\ldots,d\}$, then for a vector $\vk a\in\R^d$ we denote \EHb{by} $\vk{a}_I=(a_i, i\in I)$ a sub-block vector of $\vk a$. Similarly, if further $J \subset \{1\ldot d\}$, for  a matrix $M=(m_{ij})_{i,j\in \{1,\ldots,d\}}\in \R^{d\times d}$ we denote \EHb{by} $M_{IJ}\ccj{=M_{I,J}}=(m_{ij})_{i\in I, j\in J}$ \EHb{the} sub-block matrix of $M$
\EHb{determined by $I$ and $J$.} Further, write $M_{II}^{-1}=(M_{II})^{-1}$ for the inverse matrix of $M_{II}$ whenever it exists.
%
The next lemma stated in  \cite{HA2005} (see also  \cite{ENJH02}) is important for several definitions in the sequel.

\BEL \label{AL}
Let
$M \in \R^{d \times d},d\ge 2$ be a positive definite  matrix. 
If  $\vk{b}\in  \R ^d \setminus (-\infty, 0]^d $, then the quadratic programming problem
$$ P_M(\vk{b}): \text{minimise $ \vk{x}^\top M^{-1} \vk{x} $ under the linear constraint } \vk{x} \ge \vk{b} $$
has a unique solution $\widetilde{\vk{b}}$ and there exists a unique non-empty
index set $I\subseteq \{1, \ldots, d\}$ so that
\begin{eqnarray}  \label{eq:IJi}
\widetilde{\vk{b}}_{I}&=&
\vk{b}_{I}\not= \vk 0_I, \quad M_{II}^{-1} \vk{b}_{I}>\vk{0}_I,\\
\text{and if}&&  I^c =\{ 1, \ldots, d\} \setminus I \not=
\emptyset, \text{ then }
\widetilde {\vk{b}}_{I^c}= - ((M^{-1} )_{I^cI^c})^{-1} (M^{-1})_{I^cI} \vk b_I=  M_{I^cI}M_{II}^{-1} \vk{b}_{I}\ge \vk{b}_{I^c}.
\label{eq:hii}
\end{eqnarray}
Furthermore,
\BQN
\min_{\vk{x} \ge  \vk{b}}
\vk{x}^\top M^{-1}\vk{x}&=& \widetilde{\vk{b}}^\top M^{-1} \widetilde {\vk{b}}   =
\vk{b}_{I}^\top M_{II}^{-1}\vk{b}_{I}>0, \\
 \label{eq:new}
\vk{x}^\top M^{-1} \widetilde{ \vk{b}}&=& \vk{x}_I^\top M_{II}^{-1} \widetilde {\vk{b}}_I=
\vk{x}_I^\top M_{II}^{-1}\vk{b}_I,  \quad
\vk{x}\in \R^d.
\EQN
If $\vk{b}= b \vk{1}, b\in (0,\infty)$, then $ 2 \le \sharp\{i: i\in I\}  \le d$.
\EEL
\EHb{Hereafter}, the unique index set $I$ that defines the solution of the quadratic programming problems in question
will be referred to as the {\it essential index set}.

For any fixed $t$, let  $I(t)\subseteq \{1\ldot d\}$ be the  essential index set of the quadratic programming problem $P_{\Sigma}(\b(t))$ where
$$\vk{b}(t)=\vk{\alpha}+t\vk{\mu},\quad  t\ge 0$$
and set
$$I(t)^c = \{1\ldot d\}\setminus I(t).$$
Next, we analyze the function $g$ defined in \eqref{fgt}. \EHc{Let us briefly mention the following standard notation
for two given functions $f(\cdot)$ and $h(\cdot)$. We write
$ f(x)=h(x)(1+o(1))$  or simply $f(x)\sim h(x)$, if  $ \lim_{x \to a}  {f(x)}/{h(x)} = 1$ \ccj{ ($a\in\R\cup\{\IF\}$)}. Further, write $ f(x) = o(h(x))$, if $ \lim_{x \to a}  {f(x)}/{h(x)} = 0$.}

\BEL \label{lem:gc1}
We have   $g\in C^1(0,\IF)$. Furthermore,  $g$ is convex and it achieves its unique minimum at
\BQN\label{eq:t0}
t_{0}=   \sqrt{\frac{\vk{\alpha}_{I}^\top \Sigma^{-1}_{II}  \vk{\alpha}_{I}   }{\vk{\mu}_{I}^\top \Sigma^{-1}_{II}  \vk{\mu}_{I}  }
	}>0,
\EQN
which is given by
\BQN \label{eq:intr1}
g(t_0)&=&\inf_{t>0}\frac{1}{t} \inf_{\vk{v} \ge \vk{\alpha}+\vk{\mu} t}  \vk{v}^\top \Sigma^{-1}  \vk{v}=\frac{1}{t_0}  \vk{b}^\top_{I} \Sigma_{II}^{-1}  \vk{b}_{I},
\EQN
with
$$\vk{b}=\b(t_0)= \vk \alpha + t_0 \vk \mu$$
 and $I=I(t_0)$ being the essential index set corresponding to $P_\Sigma (\vk b)$. Moreover,
   \BQN
g(t_0\pm t) = g(t_0)+ \frac{g^{''}(t_0\pm)}{2} t^2\ooo, \ \ \ \ t\downarrow 0. \label{eq:gt0pm}
\EQN
\EEL
The proof of \nelem{lem:gc1} is displayed in the Appendix.

Hereafter we shall use the notation $\vk{b}=\b(t_0)$, and $I=I(t_0)$ for the essential index set of the quadratic programming problem $P_{\Sigma}(\b)$.  Furthermore,
let  $\widetilde{\vk{b}}$ be the \EHb{unique} solution
of $P_{\Sigma}(\b)$. 
If $I^c=\{1,\ldots,d\}\setminus I\neq\emptyset $, we define
\EHb{the} {\it weakly essential} and the {\it unessential index} sets by
\BQN\label{K.def}
K=\{j\in I^c: \widetilde{\vk{b}}_j=\Sigma_{jI}\Sigma_{II}^{-1}\vk{b}_I=\vk{b}_j\},
 \quad
\text{  and } J=\{j\in I^c: \widetilde{\vk{b}}_j=\Sigma_{jI}\Sigma_{II}^{-1}\vk{b}_I>\vk{b}_j\},
\EQN
respectively.
Set for $t>0$
\COM{
From the above discussion we can assume that the optimal sets corresponding to the quadratic programming problem $P_\Sigma(\vk b)$, with $\vk b= \vk \alpha+\vk\mu t_0$, are  $I=I(t_0)$, $K=K(t_0) $ and $J=J(t_0)$ such that
\BQNY
&&\tilb_I= \b_{I}>\vk 0_{I}\\
&&\tilb_K=\Sigma_{KI}\Sigma_{II}^{-1} \vk b_{I} = \b_{K}\\
&&\tilb_J=\Sigma_{JI}\Sigma_{II}^{-1} \vk b_{I} > \b_{J}\\
&&\Sigma_{II}^{-1}  \b_{I}>\vk 0_{I}.
\EQNY
}
\BQNY
g_I(t)=\frac{1}{t}    \vk{\alpha} ^\top_{I}  \Sigma^{-1}_{II}   \vk{\alpha} _{I} +
2 \vk{\alpha} ^\top_{I}  \Sigma^{-1}_{II}   \vk{\mu}  _{I}+
\vk{\mu}^\top_{I}  \Sigma^{-1}_{II}   \vk{\mu} _{I} t.
\EQNY
Clearly, by  \nelem{AL} we have $g(t_0)=g_I(t_0)$.
Furthermore, we have
\BQN\label{eq:gt0t1}
g_I(t_0+t)=g_I(t_0)+\frac{g_I^{''}(t_0)}{2}t^2(1+o(1)),\ \ t\to 0,
\EQN
with
$$
g_I^{''}(t_0)=2 t_0^{-3} ( \vk{\alpha}_{I}^\top \Sigma^{-1}_{II}  \vk{\alpha}_{I}).
$$
For notational simplicity we shall set below
\BQN\label{GTT}
 \gt= \inf_{t\ge 0} g(t)= g(t_0)= g_I(t_0), \quad \ggt= g_I^{''}(t_0).
 \EQN

\section{Main Results}\label{Sec:main}

Let for \EHc{the} \EHb{non-empty index set $K$ defined in (\ref{K.def}) $\vk{Y}_{K}\overset{d}\sim{\mathcal N}(\vk{0}_{K},{D}_{KK})$, i.e.,}  $\vk{Y}_{K}$ is  a normally distributed random vector with mean vector $\vk0_K$ and covariance matrix ${D}_{KK}$ \EHb{given by}
$${D}_{KK}=\Sigma_{KK}-\Sigma_{KI}
\Sigma_{II}^{-1}\Sigma_{IK}.$$
\EHb{We write}  $m=\sharp I:=\sharp\{i: i\in I\}\EHb{\ge 1}$  \EHb{for the number of elements of the index set $I$. Further}  define the following constant
\BQN\label{eq:HITT}
\HAS = \lim_{T\to\IF}\frac{1}{T}\HAS(T), \quad     \HAS(T)=\int_{\R^{m}}e^{\frac{1}{t_0}\vk{x}_I^\top\Sigma_{II}^{-1}\vk{b}_I}
\pk{\exists_{t\in[0,T]}(\vk{X}(t)-\vk{\mu}t)_I>\vk{x}_I}\,d\vk{x}_I,
\EQN
with respect to the essential index set $I$ \EHb{and set}
\BQNY
\ci= \frac{1}{\sqrt{(2\pi t_0)^m \abs{\Sigma_{II}}}} \int_{\R} e^{-\ggt \frac{x^2}{4}} \psi(x)\,dx,
\EQNY
where, \K1{for $x\inr$
\BQN\label{eq:psi}
\psi(x)=
\left\{\begin{array}{cc}
1, & \hbox{if } K=\emptyset\\
\pk{\vk{Y}_K>\frac{1}{\sqrt {t_0}}(\vk{\mu}_K-\Sigma_{KI}\Sigma_{II}^{-1}\vk{\mu}_I)x}, & \hbox{if } K\neq\emptyset .
\end{array}\right.
\EQN
\EHc{$\HAS$'s} are multidimensional counterparts of the celebrated {\it Pickands constants},
defined \EHc{in the 1-dimensional setup} as
$$\lim_{T\to\infty} \frac{1}{T} \int_{\R}e^x \pk{\exists_{t\in[0,T]} (\sqrt{2}W_H(t)-t^{2H})> x}\,dx,$$
where $W_H$ is a fractional Brownian motion with Hurst parameter $H\in(0,1]$; see also \cite{DHJT15} for the analog of
$\HAS$ when $\Sigma_{II}=I_d$ is the identity matrix.
We refer to \cite{PickandsA,PicandsB,Pit96,DRolski,DM} and references therein for properties and extensions of the
notion of classical Pickands constants.}

The next theorem constitutes our principal result. Its proof \K1{is} demonstrated in Section \ref{proofmain}.

\BT\label{Thm1} \label{th:main1}
Let $\vk \alpha, \vk \mu$ satisfy \eqref{alphamu} and let $\gt, \ggt$ be given by \eqref{GTT}. We have as $u\to \IF$
\BQN\label{eq:PuAsym}
P(u)  \sim  \ci\HAS u^{\frac{1- m}{2}}e^{- \frac{\gt}{2} u},
\EQN
where
\BQN\label{eq:HIbounds}
	  0< \frac{\EH{t_0^{m-1}} \vk{\mu}_I^\top \SI_{II}^{-1}\vk{b}_I}{16 \prod_{i\in I} (  \Sigma^{-1}_{II}\vk{b}_I )_i}\le \H_I <\IF.
\EQN
\ET
\begin{remark}
In the case that $K=\emptyset$, direct calculations show that \eqref{eq:PuAsym} holds with
\BQNY
\ci= \frac{2^{1-m/2}\pi^{(1-m)/2}}{\sqrt{  t_0^m \ggt \abs{\Sigma_{II}}}}\EHb{>0}.
\EQNY
\end{remark}

Using the same technique as in the proof of Theorem \ref{th:main1}, we can derive
the approximation of the conditional passage time $\tau_u \lvert  \tau_u<\IF$.
 	\BT\label{KorrRT}
	Let $\tau_u$ be defined  in \eqref{eq:RT} and $\psi$ be defined  in  \eqref{eq:psi}.    \ccj{Under the assumptions of \netheo{Thm1}} for any $s \in \R$ \EHb{we have}
	\BQNY
	\lim_{u\to\IF}\pk{\frac{\tau_u -t_0 u} {\sqrt{ 2u /\ggt }} \le s  \Bigl \lvert \tau_u<\IF}=\frac{\int_{-\IF}^s e^{- \frac{x^2}{2}} \psi(\sqrt{2/\ggt }x)\,dx  }{ \int_{-\IF}^\IF e^{-\frac{x^2}{2}} \psi(\sqrt{2/\ggt }x)\,dx }.
	\EQNY
	\ET

\begin{remark}
\EHb{If} $K= \emptyset$, then  \K1{by (\ref{eq:psi})}
$\frac{\tau_u -t_0 u} {\sqrt{ 2 u/\ggt} }   \Bigl \lvert \tau_u<\IF  $ is \K1{asymptotically, as $u\to\infty$,} approximated by a standard normal random variable.
\end{remark}


In the rest of this section \K1{we discuss some interesting special cases and examples.}

\subsection{Independent \EHb{components}}
\def\NN{d}
\def\mcalA{\mathcal{A}}
\def\mcalB{\mathcal{B}}

Let $\Sigma$ be the  $d\times d$ identity matrix. We focus on the case where \cc{$\vk \alpha>\vk 0$} and
$$\mu_i>0, \quad 1 \le i \le n, \quad  \mu_{j} \le 0, \quad n+1 \le j \le d $$
\EHb{is valid}  for some \EHb{positive integer} $n < d$.  The result for the easier case $n=d$ will also be included.
\EHb{Under the above assumptions} 
$$g(t)=\frac{1}{t}\inf_{\vk{v}\ge \vk{\alpha}+t\vk{\mu} } \vk{v}^\top \vk{v},  \qquad t>0.$$
Before we state the result we need \EH{to introduce} some notation.
By rearranging indexes we can have the following order of \K1{constants}
\begin{equation} \label{eq:mualpAppendix}
\frac{\abs{\mu_{n+1}}}{\alpha_{n+1}}\le \frac{\abs{\mu_{n+2}}}{\alpha_{n+2}}\le \cdots \le \frac{\abs{\mu_{ d }}}{\alpha_{ d}}.
\end{equation}
\EHc{Next,} define $d=k_1>\cdots>k_l>k_{l+1}=n$ for which
\begin{eqnarray*}
\frac{\abs{\mu_{d}}}{\alpha_{d}}&=&\frac{\abs{\mu_{j}}}{\alpha_{j}}; \qquad
k_2<j\le k_1=d\\
\vdots &&\vdots\\
\frac{\abs{\mu_{k_l}}}{\alpha_{k_l}}&=&\frac{\abs{\mu_{j}}}{\alpha_{j}}; \qquad
n+1\le j\le k_{l}
\end{eqnarray*}
and
$$
\frac{\abs{\mu_{k_{l}}}}{\alpha_{k_{l}}}
< \frac{\abs{\mu_{k_{l-1}}}}{\alpha_{k_{l-1}}}< \cdots < \frac{\abs{\mu_{ k_1 }}}{\alpha_{ k_1 }}
$$
implying that
$0=t_0^{'}<t_1^{'}<\ldots<t_l^{'}<t_{l+1}^{'}=\IF,$
where
$  t_i^{'}=\frac{\alpha_{k_{i}}}{\abs{\mu_{ k_{i} }}}, i=1,\ldots,l$
are consecutive change of dimension  instants.
Precisely, \EHb{for the quadratic programming problem in question},  constancy segments are 
$U_i=[t^{'}_{i-1},t^{'}_{i}), i=1,\ldots,l+1,$
and for $t\in U_i$ we have $I(t)=I_i=\{1,\ldots,k_i\}$ since
$$\vk{\alpha}_{I_i}+t\vk{\mu}_{I_i}> \vk{0}_{I_i}, \quad
\vk{\alpha}_{I_i^c}+t\vk{\mu}_{I_i^c}\le \vk{0}_{I_i^c}, \quad t\in U_i.$$
Define for $i=1,\ldots,l+1$ the following  auxiliary functions
\begin{eqnarray*}
g_{I_i}(t) =\frac{1}{t}\vk{\alpha}_{I_i}^\top \vk{\alpha}_{I_i}+2\vk{\alpha}_{I_i}^\top \vk{\mu}_{I_i}
+t\vk{\mu}_{I_i}^\top\vk{\mu}_{I_i}, \quad t>0
\end{eqnarray*}
and remark that $g(t)=g_{I_i}(t)$ for $t\in U_i$. Clearly,
$g_{I_i}(t), t>0$
achieves its \EHc{global} minimum at
$$
t_0(i)=\sqrt{\frac{\vk{\alpha}_{I_i}^\top \vk{\alpha}_{I_i}}{\vk{\mu}_{I_i}^\top \vk{\mu}_{I_i}}}>0.
$$
Set below
$$\EH{p:}=\min\{i=1,\ldots,l +1:\ t_{i-1}^{'}\le t_0(i)< t_{i}^{'}\}, \quad t_0=t_0(\EH{p})=\sqrt{\frac{\vk{\alpha}_{I_{\EH{p}}}^\top \vk{\alpha}_{I_{\EH{p}}}}{\vk{\mu}_{I_{\EH{p}}}^\top \vk{\mu}_{I_{\EH{p}}}}} \EHc{> 0}.$$
With the same arguments
as at the end of the proof of \nelem{lem:gc1} \EH{it follows that} $g$  achieves its minimum at $t_0$. Then
with the notation of \netheo{Thm1} we have
\BQNY
g_I(t_0)=\frac{1}{t_0}\sum_{j=1}^{k_{\EH{p}}}(\alpha_j+\mu_j t_0)^2,\ \ \  g_I^{''}(t_0) =2 t_0^{-3}\sum_{j=1}^{k_{\EH{p}}}\alpha_j ^2, \ \ \ I= I_{\EH{p}},\ \ \ m=k_{\EH{p}}.
\EQNY
Moreover, if $ t_{\EH{p}-1}^{'}< t_0(\EH{p})< t_{\EH{p}}^{'}$, then
$
K=\emptyset,
$
and  if $ t_{\EH{p}-1}^{'}=t_0(\EH{p})$, then
$
K=\{k_{\EH{p}}+1,\ldots,  k_{\EH{p}-1}\}.
$
\cc{Set}  
\BQNY
 \H_{I_{\EH{p}}}=\lim_{T\to\IF}\frac{1}{T}\int_{\R^{k_{\EH{p}}}}e^{\sum_{j=1}^{k_{\EH{p}}}(\alpha_j t_0^{-1}+\mu_j ) x_j}
\pk{\exists_{t\in[0,T]}(\vk{B}(t)-\vk{\mu}t)_{I_{p}}>\vk{x}_{I_{\EH{p}}}}\,d\vk{x}_{I_{\EH{p}}}.
\EQNY
\EHc{\Lc{We} define $\Psi(x)=1- \Phi(x), x\in \R$ with
$\Phi$ the distribution function of an $\ccj{\mathcal N}(0,1)$ random variable. Below we shall
put $\prod_{i\in \mcalB} (\cdots )=1$ for $\mcalB$ empty.}

\EHb{We reformulate next our main findings for this particular \K1{case.}}
\BK\label{cor:indep} 

(i). If $n<d$, then \cc{as $u\to\IF$}
\BQN\label{eqMth}
\pk{\exists_{t\ge 0}  \vk{B}(t)-  \vk{\mu}t> \vk{\alpha} u } \sim \frac{\H_{I_{\EH{p}}}}{\sqrt{(2\pi t_0)^{k_{\EH{p}}}  }} u^{\frac{1-k_{\EH{p}}}{2}}e^{-\frac{1}{2t_0}\sum_{j=1}^{k_{\EH{p}}}(\alpha_j+\mu_j t_0)^2 u}\int_{\R} e^{-\frac{ \sum_{j=1}^{k_{\EH{p}}}\alpha_j^2}{2t_0^{3}}x^2} \prod_{i\in K} \Psi\LT(\frac{\mu_i }{\sqrt t_0}x\RT)\,dx,
\EQN
\cc{and for any $s\in\R$
	\BQN \label{eqMth2}
\lim_{u\to\IF}	\pk{\frac{\tau_u -t_0 u} {\sqrt{  t_0^{3 } (\sum_{j=1}^{k_{\EH{p}}}\alpha_j^2)^{-1 }\EH{ u   }} } \le s  \Bigl \lvert \tau_u<\IF}= \frac{\int_{-\IF}^s e^{- \frac{x^2}{2}} \prod_{i\in K} \Psi\LT( \mu_i   t_0  (\sum_{j=1}^{k_{\EH{p}}}\alpha_j^2)^{-1/2 } x\RT) \,dx  }{ \int_{-\IF}^\IF e^{-\frac{x^2}{2}}  \prod_{i\in K} \Psi\LT( \mu_i   t_0  (\sum_{j=1}^{k_{\EH{p}}}\alpha_j^2)^{-1/2 } x\RT)\,dx }.
	\EQN }
(ii). If $n=d$, then \eqref{eqMth} and \eqref{eqMth2} hold with $\EH{p}$ replaced by $l+1$ and $K$ replaced by $\emptyset$.

\EK

\subsection{Homogeneous \EHb{$\vk \alpha$ and $\vk \mu$}}

\EHb{\K1{Suppose that} $\vk{\alpha}= \vk 1 \alpha, \alpha>0$ and
$\vk \mu= \vk 1 \mu, \mu>0$.} 
\EH{\K1{Then} for any $t>0$}
\BQNY
g(t)=  \frac{1 }{t} \inf_{\vk{v}\ge\vk{\alpha}+\vk{\mu}t}\vk{v}^\top \Sigma^{-1}\vk{v} 
=   \frac{(\alpha+\mu t)^2}{t}D,  \quad  D:= \inf_{\vk{v}\ge\vk{1} } \vk{v}^\top\Sigma^{-1}\vk{v}.
\EQNY
Let $I$ be the essential index set of the \EH{quadratic programming} problem $P_\Sigma(\vk 1)$
\K1{with $m=\sharp\{i: i\in I\}$}.  \K1{If $I^c$ is non-empty, we set
$$
K=\{j\in I^c: \Sigma_{jI}\Sigma_{II}^{-1}\vk{1}_I=\vk{1}_j\}.
$$
}
Obviously, $I(t)=I, t\ge 0$.
\EH{Further}, $g_I(t)=g(t),t>0$ and
\BQNY
t_0= \alpha \mu^{-1},\ \ 
g_I(t_0)=4D \alpha \mu, \ \ 
g_I^{''}(t_0)= 2D \mu^3 \alpha^{-1}, \ \ \cc{\vk b= 2\alpha \vk 1}.
\EQNY
\BK\label{cor:homo}
We have, as $u\to\IF$,
\BQNY
P(u) \sim \frac{2^{-\sharp K } \alpha^{(1-m)/2}\mu^{(m-3)/2} }{\sqrt{(2\pi  )^{m-1} D \abs{\Sigma_{II}}}} \H_I  u^{\frac{1-m}{2}}e^{- 2D \alpha \mu u },  
\EQNY
\cc{and for any $s\in\R$
$$\lim_{u\to\IF}\pk{\frac{\tau_u - \alpha \mu^{-1} u} {\sqrt{   \alpha D^{-1} \mu^{-3}\EH{ u   }} } \le  s \Bigl \lvert \tau_u<\IF }=\Phi(s). $$ 
}

\EK

\subsection{Negatively associated components}
In this subsection we suppose that
$$\Sigma^{-1}\vk{\alpha}>\vk{0}, \quad
\Sigma^{-1}\vk{\mu}>\vk{0}.$$
A special case of interest is \K1{when} $\Sigma^{-1}$ \EH{has all elements positive}, $\vk{\alpha}>\vk{0}$
and $\vk{\mu}>\vk{0}$.
Recall that if \EH{the} covariance matrix $\Sigma$ of $\vk{Z}$ is a correlation matrix, then
\K1{the statement that}
$\Sigma^{-1}$ \EH{has all elements positive means that} it is an  $M$-matrix, i.e.,
$\Sigma=I_d-B$, where $B\ge0$ and $I_d$ is the identity matrix.
For general  covariance matrix $\Sigma$,
\K1{ with nonpositive elements out of the diagonal}, transformation $diag(\sigma_{jj}^{-1})
\Sigma diag(\sigma_{jj}^{-1})$
makes it an \cc{$M$-matrix.} 
Notice that \K1{if a Gaussian vector $\vk{Z}$
has such \EHb{a} covariance matrix, then
$\vk{Z}$ is \it{negatively associated} (for definition
and properties see \cite{Joag-Dev83}). }

In this case
$I=\{1,\ldots,d\}$, $K=J=\emptyset$, $m=d$ and
 $$g_I(t)=
\frac{1}{t}\vk{\alpha}^\top\Sigma^{-1}\vk{\alpha}+
2\vk{\alpha}^\top\Sigma^{-1}\vk{\mu}+t \vk{\mu}^\top\Sigma^{-1}\vk{\mu}.$$
Consequently,
$$t_0=\sqrt{
\frac{\vk{\alpha}^\top\Sigma^{-1}\vk{\alpha}}{ \vk{\mu}^\top\Sigma^{-1}\vk{\mu}}} \EHc{>0}, \quad 
g_I(t_0)=\frac{1}{t_0}\vk{b}^{\top}\Sigma^{-1}\vk{b},\qquad  
g_I^{''}(t_0)=\frac{2 \vk{\alpha}^\top\Sigma^{-1}\vk{\alpha}}{ t_0^3},    $$
where $\vk{b}=\vk{\alpha}+t_0\vk{\mu}$.
Hence we \EHb{arrive at} \K1{the following result.}
\BK \EHb{As}  $u\to\infty$
$$P(u)\sim \frac{(2\pi)^{(1-d)/2}\H_{\{1,\ldots,d\}}}{\sqrt{ t_0^{d-3}  (\vk{\alpha}^\top\Sigma^{-1}\vk{\alpha}) \abs{\Sigma}}}
  u ^{\frac{1-d}{2}}e^{-\EH{\frac{\vk{b}^{\top}\Sigma^{-1}\vk{b}}{2t_0}}u},$$
\cc{and for any $s\in\R$
$$\lim_{u\to\IF}\pk{\frac{\tau_u -t_0 u} {\sqrt{  t_0^{3 }  (\vk{\alpha}^\top\Sigma^{-1}\vk{\alpha})^{-1 }\EH{ u   }} }\le s   \Bigl \lvert \tau_u<\IF } =\Phi(s).$$ 
}
\EK

\subsection{Two-dimensional case}\label{ss.two-dim}
\def\tzo{t_0^{(1)}}
\def\tzt{t_0^{(2)}}
In this section we analyze some interesting scenarios of the two-dimensional case, in which we can observe  how different entries of the covariance matrix yield different scenarios of asymptotic behaviour.
 Proofs will be postponed  to Section \ref{ss.anal.2-dim},  after presenting required results 
\K1{on} a quadratic programming problem.
For simplicity, we shall assume that
\K1{
$$\Sigma=
\left(
   \begin{array}{cc}
     1 & \rho \\
     \rho & 1 \\
   \end{array}
 \right), \quad \rho \in (-1,1)$$
and
$\mu_i=1, i=1,2,$ $\alpha_1>\alpha_2>0$.
}
\COM{Let
$
\Sigma=\left(
         \begin{array}{cc}
           1 & \rho \\
           \rho & 1 \\
         \end{array}
       \right)
,$  with inverse matrix
\BQNY
\Sigma^{-1}=\frac{1}{1-\rho^2}\left(
         \begin{array}{cc}
           1 & -\rho \\
           -\rho & 1 \\
         \end{array}
       \right).
\EQNY
}

We present next the asymptotics of \eqref{eq:Pu} for the 2-dimensional model.
\BK\label{cor:TD} (i). If $-1<\rho<\frac{\alpha_1+\alpha_2}{2\alpha_1}$, then as $u\to\IF$
\BQNY
P(u) \sim \frac{\H_{\{1,2\}} }{\sqrt{t_0^2\pi (1-\rho^2) \ggt}}   u^{-\frac{1}{2}}e^{- \frac{\gt}{2} u }, 
\EQNY
with
\BQNY
&&\ccj{t_0=\sqrt{\frac{\alpha_1^2+\alpha_2^2-2\alpha_1\alpha_2\rho}{2(1-\rho)}} \EHc{>0},\quad  \gt=\frac{2}{1+\rho}(\alpha_1+\alpha_2+2t_0 ),\quad  \ggt=2 t_0 ^{-3}\frac{\alpha_1^2+\alpha_2^2-2\alpha_1\alpha_2\rho}{1-\rho^2},}\\
&&\H_{\{1,2\}}=\lim_{T\to\IF}\frac{1}{T} \int_{\R^{2}}e^{\LT(\frac{\alpha_1-\rho\alpha_2}{(1-\rho^2) t_0}+\frac{1}{1+\rho}\RT)x_1 +\LT(\frac{\alpha_2-\rho\alpha_1}{(1-\rho^2)t_0}+\frac{1}{1+\rho}\RT)x_2 }
\pk{\exists_{t\in[0,T]}(\vk{X}(t)- t\kk{(1,1)^T}) >\vk{x} }\,d\vk{x}.
\EQNY
Furthermore, \cc{ for any $s\in\R$
$$\lim_{u\to\IF}\pk{\frac{\tau_u - t_0  u} {\sqrt{ 2u  / \ggt } } \le s  \Bigl \lvert \tau_u<\IF }=\Phi(s). $$ 
}

(ii). If $ \rho=\frac{\alpha_1+\alpha_2}{2\alpha_1}$, then as $u\to\IF$
\BQNY
P(u) \sim \frac{1 }{\sqrt{ 2\pi \alpha_1}}  \int_{\R}e^{-\frac{1}{2\alpha_1}x^2}\Psi\LT(\frac{1-\rho}{\sqrt{\alpha_1}}x\RT)\, dx e^{-2\alpha_1 u},
\EQNY
\cc{and for any $s\in\R$
	\BQNY
\lim_{u\to\IF}	\pk{\frac{\tau_u -\alpha_1 u} {\sqrt{ \alpha_1 \EH{ u   }} } \le s  \Bigl \lvert \tau_u<\IF}= \frac{\int_{-\IF}^s e^{- \frac{x^2}{2}}  \Psi\LT(  (1-\rho) x\RT) \,dx  }{ \int_{-\IF}^\IF e^{-\frac{x^2}{2}}  \Psi\LT(   (1-\rho) x\RT)\,dx }.
	\EQNY }
(iii). If $\frac{\alpha_1+\alpha_2}{2\alpha_1}<\rho<1$, then as $u\to\IF$
\BQNY
P(u) \sim e^{-2\alpha_1 u},
\EQNY
\cc{and for any $s\in\R$
$$\lim_{u\to\IF}\pk{ \alpha_1^{-1/2}    (\tau_u -  \alpha_1 u  )/\sqrt{u} \le s \Bigl \lvert \tau_u < \IF}=\Phi(s).$$} 
\EK

\begin{remark}
\EHb{According to  our findings, in both} (ii) and (iii) \EHc{above} \EH{we should  also have } the following constant
$$\H_{\{1\}}=\lim_{T\to\IF}\frac{1}{T} \int_{\R}e^{2 x }
\pk{\exists_{t\in[0,T]}(\vk{X}(t)- t) _1>{x} }\,d{x}.
$$  
However, a simple comparison with the known Pickands constants for the standard Brownian \EHb{motion}, \EHc{i.e.,}
$$\lim_{T\to\IF}\frac{1}{T}  \int_{\R}e^x\pk{\sup_{t\in[0,T]}(\sqrt{2}B_1(t)-t)>x}\,dx=1$$
yields $\H_{\{1\}}=1$.
\end{remark}

We conclude \EHc{this section} with \EHc{some} observations.

  It is possible to have similar asymptotics of $P(u)$ for $d\ge 2$ as in the 1-dimensional case. For instance \ccj{ in the above the 2-dimensional setup,}
\COM{if $\alpha_1> \alpha_2 > 0$ and
$
\Sigma=\left(
\begin{array}{cc}
1 & \rho\\
\rho & 1 \\
\end{array}
\right)
,$ then }
for \cc{$\rho\in ((\alpha_1+ \alpha_2)/(2 \alpha_1),1)$} we have  
\BQN \label{bauarb}
 \pk{\exists_{t\ge 0} \cap_{i=1}^2 \{(X_i(t) - t) > \alpha_i u\}} \sim
\ci \pk{\exists_{t\ge 0} (X_1(t) - t) > \alpha_1 u},\ \ u\to\IF,
 \EQN
with  $\ci=1$. \EHc{Consequently,}  only the first component of $\vk X(t), t\ge 0$ is controlling the asymptotics of $P(u)$. This case will be referred to  as the {\it loss of dimensions phenomena}.

  There are other cases of  {loss of dimensions phenomena}, where some components other than those with indexes in $I$ still play a role in the asymptotics of $P(u)$, but only up to some constants. For instance, referring again to the  2-dimensional case presented in \nekorr{cor:TD} we have for $\rho=(\alpha_1+ \alpha_2)/(2 \alpha_1)$ that
  \eqref{bauarb} holds, with $\ci$ taking the information of the second component and given by
  $$ \ci=   \frac{1 }{\sqrt{ 2\pi \alpha_1}}\int_{\R}e^{-\frac{1}{2\alpha_1}x^2}\Psi\LT(\frac{1-\rho}{\sqrt{\alpha_1}}x\RT)\, dx.
  $$
There are several technical issues related to the loss of dimensions as it will be explained in our proofs below.




\section{Proof\EHb{s} of Main Results} \label{proofmain} 
In this section we first  present the proof of Theorem \ref{Thm1}. \EHb{In order to convey the main ideas and to reduce complexity}, we shall divide the 
proof into several steps and \EHc{then we complete the proof  by putting all the arguments together.}

By the self-similarity of Brownian motion, for any $u$ positive we have
\BQNY
P(u)&=&\pk{\exists_{t\ge0} \ \  \vk{X}(t)-\vk{\mu} t>\vk{\alpha} u  } = \pk{\exists_{t\ge0} \ \  \vk{X}(t)>\sqrt u (\vk{\alpha}+\vk{\mu} t)}.
\EQNY
We have thus the following sandwich bounds
\BQN\label{eq:pr}
p(u)\le P(u)\le  p(u)+ r(u),
\EQN
where
\BQNY
p(u):=\pk{\exists_{t\in \Del_u}    \vk{X}(t)>\sqrt u (\vk{\alpha}+\vk{\mu} t)     }, \ \ r(u):=\pk{ \exists _{t \in \widetilde\Del_u } \vk{X}(t)> \sqrt u (\vk{\alpha}+\vk\mu t) },
\EQNY
with \EHb{(recall the definition of $t_0$ in \eqref{eq:t0})}
$$
\Del_u=\left[t_0-\frac{\ln(u)}{\sqrt{u}}, t_0+\frac{\ln(u)}{\sqrt{u}} \right],\quad
\widetilde\Del_u=\left[0,t_0-\frac{\ln (u)}{\sqrt u} \right]\cup \left[t_0+\frac{\ln (u)}{\sqrt u},\IF \right).$$

\subsection{Analysis of $r(u)$}\label{chopping off}


This step is concerned with sharp upper bounds for $r(u)$ when $u$ is large.

\BEL \label{PitI}
 For all large $u$  we have
\BQN \label{eq:t0p}
\pk{ \exists _{t \in [t_0+\frac{\ln (u)}{\sqrt u},\IF)} \vk{X}(t)> \sqrt u (\vk{\alpha}+\vk\mu t) } \le \EHb{C}    e^{- \frac{u }{2 } g_I(t_0) -  \LT(\frac{ g^{''}(t_0+) }{2}-\vn\RT) (\ln(u))^2},
\EQN
and
\BQN \label{eq:t0m}
\pk{ \exists_{ t \in [0,t_0-\frac{\ln (u)}{\sqrt u}]} \vk{X}(t)> \sqrt u (\vk{\alpha}+\vk\mu t) } \le \EHb{C} e^{- \frac{u }{2 } g_I(t_0) -  \LT(\frac{ g^{''}(t_0-) }{2}-\vn\RT) (\ln(u))^2}
\EQN
\EHb{are valid} for some constant  $ C>0$ \ccj{and some sufficiently small $ \vn>0$} which do not depend on $u$.
\EEL


{\bf Proof:}
We only present the proof of \eqref{eq:t0p} \EHb{since}  the proof of \eqref{eq:t0m} follows with similar arguments.
\LL{First note that for any $D\subset \R_+$ and any $u$ positive}
\BQNY  
 &&\pk{ \exists _{t \in D} \vk{X}(t)> \sqrt u (\vk{\alpha}+\vk\mu t) }  \\
&&\le \pk{ \exists _{t \in D}(\vk{X}(t))_{I(t)}> \sqrt u (\vk{\alpha}+\vk\mu t)_{I(t)} }  \\
&& \le \pk{ \exists _{t \in D} (\Sigma_{I(t)I(t)}^{-1} (\vk{\alpha}+\vk\mu t)_{I(t)})^\top(\vk{X}(t))_{I(t)}> \sqrt u (\Sigma_{I(t)I(t)}^{-1} (\vk{\alpha}+\vk\mu t)_{I(t)})^\top(\vk{\alpha}+\vk\mu t)_{I(t)} }\\
&&=\pk{ \exists _{t \in D} Y_{I(t)}(t)> \sqrt u  },
\EQNY
\COM{
It follows that if $\pi$ is chosen sufficiently small then $I(t)=I(t_0)=:I$ for any $t \in [t_0-\pi,t_0-\frac{\ln (u)}{\sqrt u}]$, and
$I(t)= :I_+$ for any $t \in [ t_0+\frac{\ln (u)}{\sqrt u},t_0+\pi]$. Then
\BQNY
&&\pk{ \exists t \in [t_0-\pi,t_0-\frac{\ln (u)}{\sqrt u}]: \vk{X}(t)> \sqrt u (\vk{\alpha}+\vk\mu t) }\\
&&\le \pk{ \exists t \in [t_0-\pi,t_0-\frac{\ln (u)}{\sqrt u}]: (\vk{X}(t))_I> \sqrt u (\vk{\alpha}+\vk\mu t)_I }\\
&&\le \pk{ \exists t \in [t_0-\pi,t_0-\frac{\ln (u)}{\sqrt u}]:  (\Sigma_{II}^{-1} (\vk{\alpha}+\vk\mu t)_I)^\top (\vk{X}(t))_I> \sqrt u  (\Sigma_{II}^{-1} (\vk{\alpha}+\vk\mu t)_I)^\top(\vk{\alpha}+\vk\mu t)_I }\\
&&\le \pk{ \exists t \in [t_0-\pi,t_0-\frac{\ln (u)}{\sqrt u}]: Y_I(t)> \sqrt u  }
\EQNY
}
where we used the fact that $ \Sigma_{I(t)I(t)}^{-1} (\vk{\alpha}+\vk\mu t)_{I(t)}>\vk 0_{I(t)}$ for all $t \ge0$, and
$$
Y_{I(t)}(t)= \frac{(\Sigma_{I(t)I(t)}^{-1} (\vk{\alpha}+\vk\mu t)_{I(t)})^\top (\vk{X}(t))_{I(t)}}{(\vk{\alpha}+\vk\mu t)_{I(t)} ^\top \Sigma_{{I(t)}{I(t)}}^{-1}(\vk{\alpha}+\vk\mu t)_{I(t)} },\ \ t\ge0.
$$
By the property of Brownian motion, we have \EH{almost surely}
$$
\lim_{t\to\IF} Y_{I(t)}(t) =0
$$
\EH{inplying that $Y_{I(t)}$ has bounded sample paths on [$a, \IF)$ for any $a>0$. Since further}
\BQNY
\Var\LT(Y_{I(t)}(t)\RT)=\frac{1}{g(t)}, \ \ t\ge 0
\EQNY
by the Borell-TIS inequality (see e.g., \cite{GennaBorell, AdlerTaylor, AZI,LifBook})  for any small $\theta>0$
\BQN\label{eq:Borell}
 \pk{ \exists _{t \in [t_0+\theta,\IF)} \vk{X}(t)> \sqrt u (\vk{\alpha}+\vk\mu t) }
&\le&\pk{ \exists _{t \in [t_0+\theta,\IF)} Y_{I(t)}(t)> \sqrt u  }\nonumber\\
&\le& e^{-\frac{(\sqrt u- C_0)^2}{2}\inf_{t \in [t_0+\theta,\IF)} g(t)}
\EQN
holds for all $u$ such that
$$ \sqrt u> C_0:=\E{\sup_{t \in [t_0+\theta,\IF)}  Y_{I(t)}(t)}.$$
It follows from \nelem{lem:II} that if $\theta \EHc{>0}$ is chosen sufficiently small, then for some $I^+$
$$g(t)=\frac{1}{t}    \vk{\alpha} ^\top_{I^+}  \Sigma^{-1}_{I^+I^+}   \vk{\alpha} _{I^+} +
2 \vk{\alpha} ^\top_{I^+}  \Sigma^{-1}_{I^+I^+}   \vk{\mu}  _{I^+}+
\vk{\mu}^\top_{I^+}  \Sigma^{-1}_{I^+I^+}   \vk{\mu} _{I^+} t
$$
for all $t\in (t_0,t_0+\theta)$. 
Furthermore,
$$
\E{(Y_{I^+}(t)-Y_{I^+}(s))^2}\le C_1 \abs{t-s}
$$
holds for all $s, t\in [t_0+\frac{\ln (u)}{\sqrt u},t_0+\theta]$, with some positive constant $C_1$. Thus, it follows from Piterbarg's inequality \EHc{in \cite{KEP2015}[Lemma 5.1]} \EHc{(see also \cite{Pit96}[Theorem 8.1] and \cite{MR3493177}[Theorem 3])} that
\BQN \label{eq:Pit}
 \pk{ \exists _{t \in [ t_0+\frac{\ln (u)}{\sqrt u},t_0+\theta]}\vk{X}(t)> \sqrt u (\vk{\alpha}+\vk\mu t) } & \le&
 \pk{ \exists _{t \in [ t_0+\frac{\ln (u)}{\sqrt u},t_0+\theta]} Y_{I^+}(t)> \sqrt u  }\nonumber\\
 &\le& C_2   u  e^{- \frac{u }{2 } \inf_{t\in [ t_0+\frac{\ln (u)}{\sqrt u},t_0+\theta]}g(t)}
\EQN
holds for all $u$ large, with some positive constant $C_2$ \EHc{not depending on $u$}. Moreover, for a small chosen $\theta>0$, there exists some $\vn>0$ such that
\BQN \label{eq:Dthetau}
\inf_{t\in[ t_0+\frac{\ln (u)}{\sqrt u},t_0+\theta]}g(t)\ge g(t_0)+\LT(\frac{g^{''}(t_0+)}{2}-\vn \RT)\frac{(\ln(u))^2}{u}
\EQN
\EHc{is valid} for all $u$ large with
$$
g^{''}(t_0+)=2 t_0^{-3} ( \vk{\alpha}_{I^+}^\top \Sigma^{-1}_{I^+I^+}  \vk{\alpha}_{I^+})>0.
$$
 Consequently, the claim in \eqref{eq:t0p} follows by \eqref{eq:Borell}, \eqref{eq:Pit}, \eqref{eq:Dthetau} and the fact that
$$
g(t_0) =g_I(t_0) < \inf_{t\in  [t_0+\theta,\IF) }g(t).
$$
\K1{Hence} the proof is complete.
\QED

\subsection{Analysis of $p(u)$}
We   investigate   the  asymptotics of $p(u)$ as $u \to\IF$.
Denote, for any fixed $T>0$ and $u>0$
\BQNY 
\Del_{j;u}=\Del_{j;u}(T)= [t_0+j Tu^{-1}, t_0+(j+1) Tu^{-1}],\ \ \ -N_u\le j\le N_u,
\EQNY
where $N_u=\lfloor T^{-1} \ln(u) \sqrt{u}\rfloor$
\EHb{(we denote by $\lfloor\cdot\rfloor$ the ceiling function)}.
By Bonferroni's inequality we have
\BQN\label{eq:thetaT}
p_1(u)\ge p(u) \ge p_2(u)-\Pi(u),
\EQN
where
\begin{equation*} 
p_1(u)= \sum_{j=-N_u-1}^{N_u}p_{j;u},\ \
p_2(u) = \sum_{j=-N_u}^{N_u-1}p_{j;u},\ \
\Pi(u) = \sum_{ -N_u\le j< l\le  N_u}p_{i,j;u},
\end{equation*}
with
\begin{equation*} 
p_{j;u}=\pk{\exists_{t\in\Delta_{j;u}} \vk{X}(t)>
\sqrt{u}(\vk{\alpha}+t\vk{\mu})}
\end{equation*}
and
$$
p_{i,j;u}=\pk{\exists_{t\in \Del_{i;u}}  \vk{X}(t)>\sqrt u (\vk{\alpha}+\vk{\mu} t),
 \ \exists_{t\in \Del_{j;u}} \vk{X}(t)>\sqrt u (\vk{\alpha}+\vk{\mu} t)    }.
 $$

\underline{\it Analysis of the single sum.}
We shall focus on the asymptotics of $p_1(u)$, which will be easily seen to be asymptotically equivalent to $p_2(u)$ as $u\to\IF$.
\COM{ Recall that
$$
p_1(u)= \sum_{j=-N_u}^{N_u}p_{j;u},
$$
with
\BQNY
p_{j;u}=\pk{\exists_{t\in\Delta_{j;u}} \vk{X}(t)>
\sqrt{u}(\vk{\alpha}+t\vk{\mu})}.
\EQNY
}

We first present a lemma concerning the finiteness of $\EHb{\H_I}(T)$ \Jc{defined in \eqref{eq:HITT}}, the constant that will appear in the asymptotics of $p_1(u)$.
\BEL \label{lem:HTT}
For any $T>0$ we have that $\H_I(T)<\IF$.
\EEL

{\bf Proof:} The claim follows if we can show that for any $\vk a_I>\vk 0_I$ and any $T>0$ we have
\BQN\label{eq:HTF}
\int_{\R^m}e^{ \x_I^\top  \vk a_I}
\pk{  \exists_{t\in[0,T]}   \xmu  _I> \vk{x}_I}\, d\vk{x}_I<\IF.
\EQN
\EHb{Clearly,}  it is sufficient to prove that
\BQNY
\int_{\abs{\vk x_I}>L\vk 1_I}e^{ \x_I^\top  \vk a_I}
\pk{  \exists_{t\in[0,T]}   \xmu  _I> \vk{x}_I}\, d\vk{x}_I<\IF
\EQNY
holds for some large $L$. Obviously, the above integral is the sum of a finite number of integrals with
$\vk x_I$ restricted to certain \K1{quadrants}.
Thus, without loss of generality, we may consider only the integral over $\{\vk x_{I_1}> L\vk 1_{I_1}, \vk x_{I_2}< -L \vk 1_{I_2}\}$ with $I_1\cup I_2=I$. By Borell-TIS inequality
\begin{eqnarray*}
	\lefteqn{\pk{  \exists_{t\in[0,T]}   \xmu  _I> \vk{x}_I}\le \pk{ \exists _{t \in [0, T]}
			\xmu _{I_1} >   \vk{x}_{I_1}} }\\
	&&\le \pk{\ccj{\exists_{t\in[0,T]}}\Sigma_{i\in I_1}(X_i(t)-\mu_it)\ge \Sigma_{i\in I_1}x_i}
	\le \exp\LT(-\mathcal{C}_1 (\sum_{i\in I_1}x_i -\mathcal{C}_2)^2\RT)
\end{eqnarray*}
holds for all $L$ large enough, with some positive constants $\mathcal{C}_1, \mathcal{C}_2$ which may depend on $T, \vk{\mu}$. Consequently, we may further write
\BQNY
&&\int_{\{\vk x_{I_1}> L\vk 1_{I_1}, \vk x_{I_2}\le -L \vk 1_{I_2}\}}e^{ \x_I^\top  \vk a_I}
\pk{  \exists_{t\in[0,T]}   \xmu  _I> \vk{x}_I}\, d\vk{x}_I\\
&&\le \int_{\vk x_{I_1}> L\vk 1_{I_1} }e^{ \x_{I_1}^\top  \vk a_{I_1}}
\exp\LT(-\mathcal{C}_1 (\sum_{i\in I_1}x_i -\mathcal{C}_2)^2\RT) \, d\vk{x}_{I_1}
\int_{ \vk x_{I_2}< -L \vk 1_{I_2}}e^{ \x_{I_2}^\top  \vk a_{I_2}}
\, d\vk{x}_{I_2}<\IF
\EQNY
establishing \Lc{ thus the  claim.} 
\QED

\BEL\label{lem:PL} We have as $u\to\infty$
\begin{eqnarray}\label{eq:PL0}
  p_1(u) \sim   p_2(u) \sim
\frac{1}{\sqrt{(2\pi t_0)^{m}|\Sigma_{II}|}  }  \frac{\H_I(T)}{T} u^{\frac{1- m}{2}}e^{-\frac{u}{2}g_I(t_0)}\int_{\R} e^{-\frac{g_I^{''}(t_0) x^2}{4}}\psi(x)\,dx,
\end{eqnarray}
where $\psi(x)$ is given in \ccj{\eqref{eq:psi}.}

\EEL

{\bf Proof:}
\def \cu {c_{j;u}}
By the independence of the increments property and the self-similarity of the Brownian motion,  we have
\BQN \label{Kd}
B_i(tu^{-1}+\cu)\EQD \sqrt{\cu} N_i+\frac{1}{\squ}B_i(t),\ t\in[0,T],\ \ i=1 \ldot d,
\EQN
   with $\cu=\cu(T)=t_0+j T/u,$  and $\vk N =(N_1 \ldot  N_d)$ \EHb{with independent
   $\mathcal{N}(0,1)$ components, being further} independent of $\vk B$.
Denote $\vk{Z}_{j;u}=  \sqrt{\cu} A \vk{N}$ with covariance matrix  $\Sigma_{j;u}=\cu\Sigma$ and \EHb{set} $$\vk{b}_{j;u}=\vk{b}_{j;u}(T)=\b(t_0+\frac{jT}{u})=\b+\frac{jT}{u} \vk{\mu}.$$
\EHb{\K1{Using}  \eqref{Kd} we obtain} 
\cc{\begin{eqnarray*}
p_{j;u}&=&
\pk{\exists_{t\in[t_0+\frac{jT}{u}, t_0+\frac{(j+1)T}{u}]}
  \vk{X}(t_0+\frac{jT}{u})+\vk{X}(t)-\vk{X}(t_0+\frac{jT}{u})>
\sqrt{u} (\vk{\alpha}+t\vk{\mu})}
\\
&=&\pk{\exists_{s\in[0,\frac{T}{u}]}\sqrt{c_{j;u}}A\vk{N}+{\vk{X}}(s)>
\sqrt{u}(\vk{\alpha}+(s+t_0+\frac{jT}{u})\vk{\mu})}\\
&=&\pk{\exists_{t\in[0,T]}\sqrt{c_{j;u}}A\vk{N}+\frac{1}{\sqrt{u}}
({\vk{X}}(t)-t\vk{\mu})>
\sqrt{u}(\vk{b}+\frac{jT}{u}\vk{\mu})}\\
&=&
\pk{\exists_{t\in[0,T]} \vk{Z}_{j;u}+\frac{1}{\sqrt{u}}({\vk{X}}(t)
-t\vk{\mu})>\sqrt{u}\vk{b}_{j;u}}.
\end{eqnarray*}
}
Since further
\BQNY
&& (\vk{Z}_{j;u})_I\  \overset{d}  \sim \   \mathcal{N}(\vk 0_I, (\Sigma_{j;u})_{II})\\
&&(\vk{Z}_{j;u})_{I^c } \mid    ((\vk{Z}_{j;u})_I=\vk w_I)\   \overset{d}\sim  \ \mathcal{N}( \Sigma_{I^c,I}  \Sigma_{II}^{-1}   \vk w_I,  (C_{jT,u})_{I^c} ),
\EQNY
with $(C_{jT,u})_{I^c}=\cu(\Sigma_{I^c,I^c} -\Sigma_{I^c,I}  \Sigma_{II}^{-1} \Sigma_{I,I^c} )$,
we have
\BQNY
p_{j;u}=\int_{\R^{m}}\phi_{(\SI_{j;u})_{II}}( \vk{w}_I) \pk{ \exists_{t \in [0, T]}
\begin{array}{ccc}
 \frac{1}{\sqrt{u}}   \xmu _I >
\squ (\vk{b}_{j;u})_I - \vk{w}_I    \\
(\vk{Z}_{j;u})_{I^c }  + \frac{1}{\sqrt{u}}   \xmu _{I^c}>
\squ (\vk{b}_{j;u} )_{I^c}  \Bigl\lvert     ((\vk{Z}_{j;u})_I=\vk w_I)
\end{array}
}\, d\vk{w}_I,
\EQNY
where
\BQNY
\phi_{(\SI_{j;u})_{II} }( \vk w_I)=\frac{1}{\sqrt{(2\pi)^{m}\abs{(\SI_{j;u})_{II} }}}\exp\LT(-\frac{1}{2}\vk w_I^\top(\Sigma_{j;u})_{II}^{-1}  \vk w_I \RT).
\EQNY
Using  a change of variable $\vk w_I=\squ (\vk{b}_{j;u})_I-\vk{x}_I/ \squ$ we obtain
\BQNY
&&p_{j;u}= \frac{u^{-m/2}}{\sqrt{(2\pi)^{m}\abs{(\SI_{j;u})_{II} }}}    \int_{\R^{m}}  \exp\LT(-\frac{1}{2} ( \squ  (\vk{b}_{j;u})_I- \vk{x}_I/\squ)^\top (\Sigma_{j;u})_{II}^{-1}    ( \squ  (\vk{b}_{j;u})_I- \vk{x}_I/\squ) \RT) \\
&&
\times \pk{ \exists_{t \in [0, T]}
\begin{array}{ccc}
   \xmu _I >
 \vk{x}_I    \\
\sqrt{ \cu }\vk{Y}_{I^c } + \Sigma_{I^c,I}  \Sigma_{II}^{-1}  ( \squ  (\vk{b}_{j;u})_I- \vk{x}_I/\squ)
+ \frac{1}{\sqrt{u}}   \xmu _{I^c}>
\squ (\vk{b}_{j;u} )_{I^c}
\end{array}
}\, d\vk{x}_I,
\EQNY
where
\BQNY
\vk{Y}_{I^c } \  \overset{d}\sim \  \mathcal{N}(\vk 0_{I^c}, D_{I^cI^c}), \ \ \ D_{I^cI^c}=\Sigma_{I^cI^c}-\Sigma_{I^cI}
\Sigma_{II}^{-1}\Sigma_{II^c}.
\EQNY
 Next, we  work out the exponent under the above integral
 \BQNY
 &&( \squ  (\vk{b}_{j;u})_I- \vk{x}_I/\squ)^\top (\Sigma_{j;u})_{II}^{-1}    ( \squ  (\vk{b}_{j;u})_I- \vk{x}_I/\squ) \\
 &&=u \frac{1}{\cu} (\vk{b}_{j;u})_I^\top\Sigma_{II}^{-1} (\vk{b}_{j;u})_I - 2 \frac{1}{\cu} \vk x_I^\top\Sigma_{II}^{-1} (\vk{b}_{j;u})_I  + \frac{1}{u \cu} \vk{x}_I^\top\Sigma_{II}^{-1} \vk{x}_I \\
 &&=  u   g_I(t_0+\frac{jT}{u})  - 2 \frac{1}{\cu} \vk x_I^\top\Sigma_{II}^{-1} (\vk{b}_{j;u})_I  + \frac{1}{u \cu} \vk{x}_I^\top\Sigma_{II}^{-1} \vk{x}_I.
 \EQNY
\COM{ 

We now work out the exponent under the integral of \eqref{eq:PL1}. Since
$$(\Sigma_{j;u})_{II}^{-1}=\frac{1}{c_{j;u}}(A^\top A)_{II}^{-1}=\frac{1}{c_{j;u}}\Sigma_{II}^{-1}$$
we have
\begin{eqnarray}
(\sqrt{u}(\vk{b}_{j;u})_I-\frac{\vk{x}_I}{\sqrt{u}})^\top
  (\Sigma_{j;u})_{II}^{-1}
  (\sqrt{u}(\vk{b}_{j;u})_I-\frac{\vk{x}_I}{\sqrt{u}})\nonumber\\
= u\frac{1}{c_{j;u}}(\vk{b}_{j;u})_I^\top\Sigma_{II}^{-1}(\vk{b}_{j;u})_I
 & -&\frac{2}{c_{j;u}}
  \vk{x}_I^\top(\Sigma_{II})^{-1}(\vk{b}_{j;u})_I
  +\frac{1}{uc_{j;u}}\vk{x}_I^{\top} \Sigma_{II}^{-1}\vk{x}_I
\nonumber\\
=ug_I(t_0+\frac{jT}{u})
 &-&\frac{2}{c_{j;u}}
  \vk{x}_I^\top(\Sigma_{II})^{-1}(\vk{b}_{j;u})_I
  +O(\frac{jT}{u}) .\label{eq:PL4} 
\end{eqnarray}
Using that
\begin{equation}\label{eq:PLtaylor}
g_I(t_0+t)=g_I(t_0)+\frac{g_I^{''}(t_0)}{2}t^2(1+o(1)),\ \ t\to 0
\end{equation}
and
\BQNY
\frac{2}{c_{j;u}}
\vk{x}_I^\top(\Sigma_{II})^{-1}(\vk{b}_{j;u})_I=
\frac{2}{t_0}
  \vk{x}_I^\top(\Sigma_{II})^{-1}(\vk{b})_I+O(\frac{jT}{u})
\EQNY
we can rewrite \eqref{eq:PL4} in the from:
\begin{equation}\label{eq:PL6} ug_I(t_0)+\frac{g^{''}(t_0)}{2}(\frac{jT}{u})^2
 -\frac{2}{t_0}
  \vk{x}_I^\top(\Sigma_{II})^{-1}(\vk{b})_I
  +O(\frac{jT}{u}).
\end{equation}
}

 Note that
\begin{eqnarray*} 
\sqrt{u}(\vk{b}_{j;u})_{I^c}&=&\sqrt{u}\vk{b}_{I^c}+\vk{\mu}_{I^c}
\frac{jT}{\sqrt{u}},\\
\sqrt{u}\Sigma_{I^c,I}\Sigma_{II}^{-1}
 (\vk{b}_{j;u})_I&=&\sqrt{u}\Sigma_{I^c,I}\Sigma_{II}^{-1}\vk{b}_I+
\frac{jT}{\sqrt{u}}\Sigma_{I^c,I}\Sigma_{II}^{-1}\vk{\mu}_I.
\end{eqnarray*}
Furthermore, denote
\begin{eqnarray*}\vk{Z}_K(t,\vk{x}_I)&=&(\vk{X}(t)-t\vk{\mu})_K-\Sigma_{KI}\Sigma_{II}^{-1}\vk{x}_I,\\
  \vk{Z}_J(t,\vk{x}_I)&=&(\vk{X}(t)-t\vk{\mu})_J-\Sigma_{JI}\Sigma_{II}^{-1}\vk{x}_I.
  \end{eqnarray*}
\EHb{For any $u$ positive we have}
\begin{eqnarray*}
\lefteqn{ \Bigl\{ \sqrt{c_{j;u}}Y_{I^c}+
\Sigma_{I^c,I}\Sigma_{II}^{-1}
(\sqrt{u}(\vk{b}_{j;u})_I-
\frac{\vk{x}_I}{\sqrt{u}})
+\frac{1}{\sqrt{u}}(\vk{X}(t)-t\vk{\mu})_{I^c}>
\sqrt{u}(\vk{b}_{j;u})_{I^c}\Bigr\}}\\
&=&
\left\{
\begin{array}{ll}
  \sqrt{c_{j;u}} Y_{K}
  +\frac{1}{\sqrt{u}}\vk{Z}_K(t,\vk{x}_I)
  >\frac{jT}{\sqrt{u}}
(\vk{\mu}_K-\Sigma_{KI}\Sigma_{II}^{-1}\vk{\mu}_I)\\
  \sqrt{c_{j;u}} Y_{J}
   +\frac{1}{\sqrt{u}}\vk{Z}_J(t,\vk{x}_I)
>\sqrt{u}(\vk{b}_J-\Sigma_{JI}\Sigma_{II}^{-1}\vk{b}_I+
(\vk{\mu}_J-\Sigma_{JI}\Sigma_{II}^{-1}\vk{\mu}_I)\frac{jT}{u})
\end{array}
\right\},
\end{eqnarray*}
where we used  $\vk{b}_K-\Sigma_{KI}\Sigma_{II}^{-1}\vk{b}_I=\vk{0}_K$.
Consequently, for the single sum we have
\BQN \label{eq:rTu}
p_1(u)&=& \frac{u^{-m/2}}{\sqrt{(2\pi)^{m}\abs{\SI_{II} }}}
  \sum_{{-N_u-1}\le j\le N_u} \frac{1}{\cu^{m/2}}\exp\LT(-\frac{1}{2} u   g_I(t_0+\frac{jT}{u})  \RT)
 \int_{\R^{m}}  f_{j;u}(T,\vk x_I)P_{j;u}(T,\vk x_I)\,d\vk{x}_I\nonumber\\
 &=:&\frac{{1}}{T}\frac{1}{\sqrt{(2\pi)^{m}\abs{\SI_{II} }}} u^{(1-m)/2}  e^{ - \frac{u   g_I(t_0)  }{2}  }
\EH{ R_T(u)},
\EQN
where
\begin{eqnarray}
f_{j;u}(T,\vk x_I)&=&  \exp\LT(\frac{1}{\cu} \vk x_I^\top\Sigma_{II}^{-1} (\vk{b}_{j;u})_I  - \frac{1}{2u \cu} \vk{x}_I^\top\Sigma_{II}^{-1} \vk{x}_I    \RT) ,\label{eq:fju}\\
P_{j;u}(T,\vk x_I) &=&\pk{
    \exists_{t\in[0,T]}
      \begin{array}{l}
        (\vk{X}(t)-t\vk{\mu})_I>\vk{x}_I\nonumber\\
        \sqrt{c_{j;u}} Y_{K}
+\frac{1}{\sqrt{u}}\vk{Z}_K(t,\vk{x}_I) >\frac{jT}{\sqrt{u}}
(\vk{\mu}_K-\Sigma_{KI}\Sigma_{II}^{-1}\vk{\mu}_I)\nonumber\\
\sqrt{c_{j;u}} Y_{J}+\frac{1}{\sqrt{u}}
\vk{Z}_J(t,\vk{x}_I)
>\sqrt{u}(\vk{b}_J-\Sigma_{JI}\Sigma_{II}^{-1}\vk{b}_I+
(\vk{\mu}_J-\Sigma_{JI}\Sigma_{II}^{-1}\vk{\mu}_I)\frac{jT}{u})
      \end{array}}.\label{eq:Pju}
\EQN
and
\BQN
R_T(u)&=& \exp\LT( \frac{u   g_I(t_0)  }{2}  \RT) \frac{T}{\sqrt{u}}\sum_{\kk{-N_u-1}\le j\le N_u}
\frac{1}{\cu^{m/2}}\exp\LT(-\frac{1}{2} u   g_I(t_0+\frac{jT}{u})  \RT)
\nonumber\\
&&\qquad\quad
\times \int_{\R^{m}}
f_{j;u}(T,\vk x_I)P_{j;u}(T,\vk x_I)\,d\vk{x}_I.
\label{eq:bPL}
\end{eqnarray}
We \EHb{shall prove} in Section \ref{TDT}  \EHb{that}
\BQN \label{rT}
\limit{u} R_T(u)&=& t_0^{-m/2}\H_I(T)
 \int_{-\infty}^\infty e^{-\frac{g_I^{''}(t_0) x^2}{4}}\psi(x)\,dx
\EQN
implying thus \eqref{eq:PL0} \EHc{(recall} that $\H_I(T)<\IF$ by \nelem{lem:HTT}).
\COM{\begin{eqnarray*}
  \EH{p_1(u)}   \sim
\frac{1}{\sqrt{(2\pi t_0)^{m}|\Sigma_{II}|}  }  \frac{\H_I(T)}{T} u^{\frac{1- m}{2}}e^{-\frac{u}{2}g_I(t_0)}\int_{-\infty}^\infty e^{-\frac{g_I^{''}(t_0) x^2}{4}}\psi(x)\,dx,
\end{eqnarray*}
and the claim follows then by
\eqref{lem:Pick-finite}, letting $T\to\IF$. }
\QED

\COM{
\BEL\label{lem:bPL}
With $ f_{j;u}(T,\vk x_I)$ and $P_{j;u}(T,\vk x_I)$ given in \eqref{eq:fju} and \eqref{eq:Pju} we have
\begin{eqnarray}\label{eq:bPL}
&&\lim_{u\to\IF}\exp\LT( \frac{u   g_I(t_0)  }{2}  \RT) \frac{T}{\sqrt{u}}\sum_{-N_u\le j\le N_u}
 \frac{1}{\cu^{m/2}}\exp\LT(-\frac{1}{2} u   g_I(t_0+\frac{jT}{u})  \RT)
     \nonumber\\
  &&\qquad\quad
 \times \int_{\R^{m}}
 f_{j;u}(T,\vk x_I)P_{j;u}(T,\vk x_I)\,d\vk{x}_I\\
 &&=\frac{1}{t_0^{m/2}}\H_I(T)
 \int_{-\infty}^\infty e^{-\frac{g_I^{''}(t_0) y^2}{4}}\psi(x)\,dx.\nonumber
\end{eqnarray}
\COM{where
\begin{equation}
\H_I(T)=\int_{\R^{m}}e^{\frac{1}{t_0}\vk{x}_I\Sigma_{II}^{-1}\vk{b}_I}
\pk{\exists_{t\in[0,T]}(\vk{X}(t)-\vk{\mu})_I>\vk{x}_I}\,d\vk{x}_I.
\end{equation}
}
\EEL
}

We \EHb{shall conclude} this section with a result which is needed to prove the sub-additivity property of $\H(T),T>0$. \EHb{In the following for any fixed $S\in \R, T>0$ we set}
$$\Del_u(S,T)=\left[t_0+Su^{-1}, t_0+(S+T)u^{-1}\right].$$
\EHb{Note in passing} that if $S=T$, then $\Del_u(S,T)=\Del_{1;u}$.
\BEL\label{lem:RPL}
  For any fixed $S\in \R, T>0$, we have as $u\to\IF$
\BQN\label{eq:TPL}
\pk{\exists_{t\in\Del_u(S,T)}  (\vk X(t))_I>\sqrt  u (\vk\alpha_I+t\vk\mu_I) }
\sim
  \frac{\pk{\vk{Y}_K>\vk{0}_K}}{\sqrt{(2\pi t_0)^{m}|\Sigma_{II}|}}\H_I(T)   u^{-\frac{m}{2}}e^{-\frac{u}{2}g_I(t_0)}
  .
  \EQN
\EEL

{\bf Proof:}  \EHc{As} in \eqref{eq:rTu} \EHc{for all $u>0$ we have}
\BQNY
&&\pk{\exists_{t\in\Del_u(S,T)}  (\vk X(t))_I>\sqrt  u (\vk\alpha_I+t\vk\mu_I) }\\
&&= \frac{u^{-m/2}}{\sqrt{(2\pi)^{m}\abs{\SI_{II} }}}
\frac{1}{(c_u(S))^{m/2}}\exp\LT(-\frac{1}{2} u   g_I(t_0+\frac{S}{u})  \RT)
 \int_{\R^{m}}  f_{u}(S,\vk x_I)P_{u}(S,T,\vk x_I)\,d\vk{x}_I,
\EQNY
where $c_u(S)=t_0+S/u,$ and with $\vk b_u(S)=\vk b +\vk \mu S/u$
\begin{eqnarray*}
f_{u}(S,\vk x_I)&=&  \exp\LT(\frac{1}{ c_u(S)} \vk x_I^\top\Sigma_{II}^{-1} (\vk{b}_{u}(S))_I  - \frac{1}{2u c_u(S)} \vk{x}_I^\top\Sigma_{II}^{-1} \vk{x}_I    \RT), \\
P_{u}(S,T,\vk x_I) &=&\pk{
    \exists_{t\in[0,T]}
      \begin{array}{l}
        (\vk{X}(t)-t\vk{\mu})_I>\vk{x}_I\nonumber\\
        \sqrt{c_{u}(S)} Y_{K}
+\frac{1}{\sqrt{u}}\vk{Z}_K(t,\vk{x}_I) >\frac{S}{\sqrt{u}}
(\vk{\mu}_K-\Sigma_{KI}\Sigma_{II}^{-1}\vk{\mu}_I)\nonumber\\
\sqrt{c_{u}(S)} Y_{J}+\frac{1}{\sqrt{u}}
\vk{Z}_J(t,\vk{x}_I)
>\sqrt{u}(\vk{b}_J-\Sigma_{JI}\Sigma_{II}^{-1}\vk{b}_I+
(\vk{\mu}_J-\Sigma_{JI}\Sigma_{II}^{-1}\vk{\mu}_I)\frac{S}{u})
      \end{array}}.
\EQNY
We adopt the same notation introduced in \eqref{eq:vnx1} and \eqref{eq:vnx2}. Next, we have the following upper bounds:
	\BQNY
		f_u(S,\vk{x}_I)  \le e^{  \frac{  \x _{I }^\top (\SI_{I I }^{-1}\b_I + \vk\vn_I^{\x_I})}{t_0+\vn(\x_I)}  },\ \ P_{u}(S,T,\vk x_I) \le \pk{
    \exists_{t\in[0,T]}
        (\vk{X}(t)-t\vk{\mu})_I>\vk{x}_I}.
		\EQNY
Furthermore,  by
\eqref{eq:HTF}
\BQNY
 \int_{\R^{m}} e^{  \frac{  \x _{I }^\top (\SI_{I I }^{-1}\b_I + \vk\vn_I^{\x_I})}{t_0+\vn(\x_I)}  } \pk{
    \exists_{t\in[0,T]}
        (\vk{X}(t)-t\vk{\mu})_I>\vk{x}_I}\,d\vk{x}_I<\IF.
\EQNY
Consequently, the claim follows from \EHb{the} dominated convergence theorem by letting $u\to\IF$, and thus the proof is complete. \QED


\underline{\it Finiteness and positivity of $\H_I$.} 
Recall that $I$ with \EHb{$m=\sharp I$ elements} is the essential index set of the quadratic programming problem $P_\Sigma(\vk b)$ where
$$\vk b=\vk{b}(t_0)=\vk\alpha+\vk \mu t_0.$$
%
\EH{We first prove the sub-additivity} of $\H_I(T), T>0$.

\BEL\label{lem:HT}
For any $S,T$ positive we have
$\HAS(S+T)\le \HAS(S)+\HAS(T)$. 
Moreover,
$$
\H_I=  \inf_{T > 0}\frac{1}{T}\H_I(T)<\IF.
$$
\EEL

{\bf Proof:}
Note that
\BQNY
 \pk{\exists_{t\in[t_0, t_0+(S+T)u^{-1}]}  (\vk X(t))_I>\sqrt  u (\vk\alpha_I+t\vk\mu_I) }
&\le &\pk{\exists_{t\in[t_0,t_0+Su^{-1}]}  (\vk X(t))_I>\sqrt  u (\vk\alpha_I+t\vk\mu_I) }\\
&&+\pk{\exists_{t\in[t_0+Su^{-1}, t_0+(S+T)u^{-1}]}  (\vk X(t))_I>\sqrt  u (\vk\alpha_I+t\vk\mu_I) }.
\EQNY
Using the result of Lemma \ref{lem:RPL} the proof of the sub-additivity follows.  The second claim follows directly from Fekete's  lemma. This completes the proof.
\QED

\BEL \label{TR}  For any $t>0$
\BQN\label{eq:const}
 \int_{\R^m}e^{\frac{\x_I^\top \SI_{II}^{-1}\vk{b}_I}{ t_0}}
 \pk{     \xmu  _I> \vk{x}_I}\, d\vk{x}_I=
 \frac{t_0^m}{ \prod_{\EH{i\in I}}(  \Sigma^{-1}_{II}\vk{b}_I )_i}>0.
\EQN
\EEL

{\bf Proof:} First note that the solution of the quadratic programming problem $P_\Sigma(\vk b)$ is such that
$$\prod_{\EH{i\in I}}(  \Sigma^{-1}_{II}\vk{b}_I )_i>0.$$
Since
$$\E{e^{\vk{s}_I^\top  (\vk{X}(t))_I}}= e^{t \vk{s}_I^\top \Sigma_{II}  \vk{s}_I/2}, \quad \EHc{\vk{s} \inr^d}, t>0$$
\EH{for any $a>0$}  we have
\BQNY
\int_{\R^m}e^{\EH{a}{\x_I^\top \SI_{II}^{-1}\vk{b}_I}}
\pk{     \xmu  _I> \vk{x}_I}\, d\vk{x}_I
&=& e^{-{\EH{a}t\vk{\mu}_I^\top  \SI_{II}^{-1}\vk{b}_I}}  \int_{\R^m} e^{\EH{a}{\vk{y}_I^\top \SI_{II}^{-1} \vk{b}_I}}\Biggl( \int_{\vk{z}_I\ge \vk{y}_I}
\phi_{ t \Sigma_{II} }(\vk{z}_I) \, d\vk{z}_I \Biggr)\, d\vk{y}_I\\
&=&\frac{a^{-m}}{ \prod_{i\in I} (  \Sigma^{-1}_{II}\vk{b}_I )_i}  e^{-{\EH{a} \vk{\mu}_I^\top  \SI_{II}^{-1}\vk{b}_I}}\int_{\R^m}  e^{\EH{a}{\vk{z}_I^\top \SI_{II}^{-1} \vk{b}_I}}  \phi_{ t \Sigma_{II} }(\vk{z}_I) \, d\vk{z}_I\\
&=& \frac{a^{-m}}{ \prod_{i\in I} (  \Sigma^{-1}_{II}\vk{b}_I )_i}  e^{-{\EH{a}t\vk{\mu}_I^\top  \SI_{II}^{-1}\vk{b}_I}
	+  a^2 t\frac{\vk{b}_I^\top  \SI_{II}^{-1}\vk{b}_I}{2}},
\EQNY
\cc{where
	$$
	\phi_{ t \Sigma_{II} }(\vk{z}_I) =\frac{1}{\sqrt{(2\pi t)^m\abs{\Sigma_{II}} }}\exp\LT(-\frac{1}{2 t} \vk{z}_I^\top \Sigma_{II}^{-1}\vk{z}_I\RT).
	$$}
\EH{In view of \eqref{eq:gt0t1} we have that $g_I'(t_0)=0$,  (recall that  $\vk{b}=\vk\alpha+\vk\mu t_0$) hence
	\BQN \label{eq:mub}
	-\frac{ \vk{\mu}_I^\top  \SI_{II}^{-1}\vk{b}_I}{t_0} + \frac{ \vk{b}_I^\top  \SI_{II}^{-1}\vk{b}_I}{2t_0^2} =0
	\EQN
	implying thus $\vk{\mu}_I^\top \SI_{II}^{-1}\vk{b}_I>0.$ Moreover, choosing
	$a=1/t_0$, where
	$$t_0=    \sqrt{\frac{\vk{\alpha}_{I}^\top \Sigma^{-1}_{II}  \vk{\alpha}_{I}   }{\vk{\mu}_{I}^\top \Sigma^{-1}_{II}  \vk{\mu}_{I}  }}> 0$$
	 establishes the claim.
}
\QED

\COM{
\EH{for ans $a>0$}
\BQNY
\int_{\R^m}e^{\EH{a}{\x_I^\top \SI_{II}^{-1}\vk{b}_I}}
\pk{     \xmu  _I> \vk{x}_I}\, d\vk{x}_I
&=& e^{-{\EH{a}t\vk{\mu}_I^\top  \SI_{II}^{-1}\vk{b}_I}}  \int_{\R^m} e^{\EH{a}{\vk{y}_I^\top \SI_{II}^{-1} \vk{b}_I}}\Biggl( \int_{\vk{z}_I\ge \vk{y}_I}
\phi_{ t \Sigma_{II} }(\vk{z}_I) \, d\vk{z}_I \Biggr)\, d\vk{y}_I\\
&=&\frac{a^m}{ \prod_{i=1}^m(  \Sigma^{-1}_{II}\vk{b}_I )_i}  e^{-{\EH{a} t\vk{\mu}_I^\top  \SI_{II}^{-1}\vk{b}_I}}\int_{\R^m}  e^{\EH{a}{\vk{z}_I^\top \SI_{II}^{-1} \vk{b}_I}}  \phi_{ t \Sigma_{II} }(\vk{z}_I) \, d\vk{z}_I\\
&=& \frac{a^m}{ \prod_{i=1}^m(  \Sigma^{-1}_{II}\vk{b}_I )_i}  e^{-{\EH{a}\vk{\mu}_I^\top  \SI_{II}^{-1}\vk{b}_I}} e^{ a^2\frac{t\vk{b}_I^\top  \SI_{II}^{-1}\vk{b}_I}{2}},
\EQNY
\cc{where
	$$
	\phi_{ t \Sigma_{II} }(\vk{z}_I) =\frac{1}{\sqrt{(2\pi t)^m\abs{\Sigma_{II}} }}\exp\LT(-\frac{1}{2 t} \vk{z}_I^\top \Sigma_{II}^{-1}\vk{z}_I\RT).
	$$}
By the properties of the quadratic programming problem $P_\Sigma(\vk b)$ we have that
$$\prod_{\EH{i\in I}}(  \Sigma^{-1}_{II}\vk{b}_I )_i>0,$$
hence the claim follows since  further \EH{$t_0=    \sqrt{\frac{\vk{\alpha}_{I}^\top \Sigma^{-1}_{II}  \vk{\alpha}_{I}   }{\vk{\mu}_{I}^\top \Sigma^{-1}_{II}  \vk{\mu}_{I}  }}> 0$.
}
	\QED
	}

\BEL \label{Lem2} 
We have
$$
	\H_I \ge \frac{\EH{t_0^{m-1}} \vk{\mu}_I^\top \SI_{II}^{-1}\vk{b}_I}{16 \prod_{i\in I} (  \Sigma^{-1}_{II}\vk{b}_I )_i} > 0.$$

\EEL
{\bf Proof:}
Suppose \K1{that} $\delta >0$ and let $n$ be any integer.  \EHb{Application of Bonferroni's inequality yields}
\begin{eqnarray*}
	\HAS(\delta n) &\ge& \int_{\R^m} e^{\frac{\x_I^\top \SI_{II}^{-1}\vk{b}_I}{ t_0}}
	\pk{\exists k\in \{1 \ldot n \}:\    ( \vk{X}(\delta k)-\vk{\mu}(\delta k))_I> \vk{x}_I}\, d\vk{x}_I\\
	&\ge&\int_{\R^m} e^{\frac{\x_I^\top \SI_{II}^{-1}\vk{b}_I}{ t_0}}
	\sum_{k=1}^n\pk{ ( \vk{X}(\delta k)-\vk{\mu}(\delta k))_I> \vk{x}_I} \, d\x _I \\
	&&\ \
	-\int_{\R^m} e^{\frac{\x_I^\top \SI_{II}^{-1}\vk{b}_I}{ t_0}}
	\sum_{k=1}^{n-1}\sum_{l=k+1}^n\pk{ ( \vk{X}(\delta k)-\vk{\mu}(\delta k))_I> \vk{x}_I,
		(\vk{X}(\delta l)-\vk{\mu}(\delta l))_I> \vk{x}_I}
	\, d\x_I  \cr
	&=:&\mathcal{I}_1- \mathcal{I}_2.
\end{eqnarray*}
By \nelem{TR} we have
$$\mathcal{I}_1=n Q, \quad Q:= \frac{t_0^m}{ \prod_{i=1}^m(  \Sigma^{-1}_{II}\vk{b}_I )_i}.$$
Next, since
\begin{eqnarray*}
	\pk{ ( \vk{X}(\delta k)-\vk{\mu}(\delta k))_I> \vk{x}_I,
		(\vk{X}(\delta l)-\vk{\mu}(\delta l))_I> \vk{x}_I}
	\le
	\pk{ \frac{1}{2}( \vk{X}(\delta k)+\vk{X}(\delta l) - \vk{\mu} (\delta k  + \delta l))_I>\x_I}
\end{eqnarray*}
by \nelem{TR} 
\BQNY
\mathcal{I}_2&\le &  \sum_{k=1}^{n-1}\sum_{l=k+1}^n
\int_{\R^m} e^{\frac{\x_I^\top \SI_{II}^{-1}\vk{b}_I}{ t_0}}
\pk{ \frac{1}{2}( \vk{X}(\delta k)+\vk{X}(\delta l) - \vk{\mu} (\delta k  + \delta l))_I>\x_I}
\, d\x_I   \\
&= &  \sum_{k=1}^{n-1}\sum_{l=k+1}^n
\int_{\R^m} e^{\frac{\x_I^\top \SI_{II}^{-1}\vk{b}_I}{ t_0}}
\pk{ \sqrt{\frac{3\delta k+\delta l}{4}}(\vk X(1))_I- \vk{\mu} _I \frac{3\delta k+\delta l}{4} >\x_I+ \frac{\vk{\mu}_I \delta (l-k) }{4} }
\, d\x_I   \\
&= &  Q \sum_{k=1}^{n-1}\sum_{l=k+1}^n e^{ -\frac{  \vk{\mu}_I^\top \SI_{II}^{-1}\vk{b}_I  \delta   }{4t_0}(l-k)  }.
\EQNY
\EH{Since by \eqref{eq:mub}  we have $\vk{\mu}_I^\top \SI_{II}^{-1}\vk{b}_I>0$, hence}
\BQNY
\mathcal{I}_2&\le  &  Q n\int_0^\IF  e^{ -\frac{  \vk{\mu}_I^\top \SI_{II}^{-1}\vk{b}_I  \delta  }{4t_0} x } \, dx
=  Q n \frac{4 t_0}{\delta   \vk{\mu}_I^\top \SI_{II}^{-1}\vk{b}_I }.
\EQNY
 By \nelem{lem:HT}
\kk{
$$ \H_I=\inf_{T > 0}\frac{1}{T}\H_I(T)
\ge\inf_{n >0}\frac{\mathcal{I}_1- \mathcal{I}_2}{\delta n}
=\frac{Q}{\delta } \left(1-  \frac{1}{\delta } \frac{4 t_0}{   \vk{\mu}_I^\top \SI_{II}^{-1}\vk{b}_I }\right).
$$
Since $\delta>0$ was arbitrary, as in \cite{DHJT15},
$$ \H_I
\ge  \max_{\delta >0}\frac{Q}{\delta } \left(1-  \frac{1}{\delta } \frac{4 t_0}{   \vk{\mu}_I^\top \SI_{II}^{-1}\vk{b}_I }\right)
\ge Q \frac{\vk{\mu}_I^\top \SI_{II}^{-1}\vk{b}_I } {16t_0}>0,$$
}
establishing the proof.  \QED

\def\wHAS{\widetilde{\HAS}}

\underline{\it Estimation of double-sum.} 
In this subsection we shall show that  as $u\to\IF$ and then $T\to\IF$
\BQN\label{eq:Pip}
\Pi(u)=o(p_1(u)).
\EQN
First, note that
\begin{eqnarray}
  \lefteqn{p_{i,j;u}=\pk{\exists_{s\in\Delta_{i;u}} \vk{X(s)}>\sqrt{u}(\vk{\alpha}+s\vk{\mu}),
    \exists_{t\in\Delta_{j;u}} \vk{X(t)}>\sqrt{u}(\vk{\alpha}+t\vk{\mu})}}\nonumber\\
&&\le
  \pk{\exists_{(s,t)\in\Delta_{i;u}\times \Delta_{j;u}} (\vk{X}(s)+\vk{X}(t))_I>\sqrt{u}(2\vk{\alpha}_I+(s+t)\vk{\mu}_I)}\nonumber\\
  &&=
    \pk{\exists_{(s,t)\in[0,T]^2}\frac{1}{2}
      (\vk{X}(t_0+\frac{iT+s}{u} )+\vk{X}(t_0+\frac{jT+t}{u}))_I>
      \sqrt{u}(\vk{\alpha}_I+(t_0+\frac{(i+j)T+s+t}{2u})\vk \mu_I)}.\label{eq:DS1}
    \end{eqnarray}
Next we rewrite for $(s,t)\in[0,T]^2$
\begin{eqnarray*}
  \lefteqn{\vk{X}(t_0+\frac{iT+s}{u})+\vk{X}(t_0+\frac{jT+t}{u})}\\
  &&=\left\{2\vk{X}(t_0+\frac{iT}{u})\right\}+
 \left \{(\vk X(t_0+\frac{iT+s}{u})-\vk{X}(t_0+\frac{iT}{u}))
  +(\vk X(t_0+\frac{(i+1)T}{u})-\vk{X}(t_0+\frac{iT}{u})) \right\}\\
  &&\qquad
  + \left\{\vk{X}(t_0+\frac{jT}{u})-\vk{X}(t_0+\frac{(i+1)T}{u}) \right\}
  +  \left \{\vk{X}(t_0+\frac{jT+t}{u})-\vk{X}(t_0+\frac{jT}{u}) \right\}.
\end{eqnarray*}
Note  that all the processes (or random variables) inside consecutive $\{\ldots\}$
are mutually independent. Consequently,
\begin{eqnarray}\label{eq:Ziju}
  \lefteqn{\frac{1}{2} \Bigl(\vk{X}(t_0+\frac{iT+s}{u})+\vk{X}(t_0+\frac{jT+t}{u})\Bigr)_I}\\
  &&\overset{d}= (\vk Z_{i,j;u})_I+\frac{1}{2\sqrt u} \Bigl(\vk{X}_1(s)+\vk{X}_1(T)+
    \vk{X}_2(t) \Bigl)_I,\quad (s,t)\in[0,T]^2,\nonumber
\end{eqnarray}
where $\vk{X}_1$ and $\vk{X}_2$ are independent copies of
$\vk{X}$, which are also independent of $\vk Z_{i,j;u}\overset{d}= \sqrt{c_{i,j;u}}A\vk{N}$, with
covariance matrix  $\Sigma_{i,j;u}=c_{i,j;u}\Sigma$, where
$c_{i,j;u}=t_0+\frac{(j+3i-1)T}{4u}$ and
$\vk N =(N_1 \ldot  N_d)$ \EHb{has independent $\mathcal{N}(0,1)$ components.}

 \cc{
Next, set
$$\vk{b}_{i,j;u}=\vk{b}\left(t_0+\frac{(i+j)T}{2u}\right).$$
It follows from
 \eqref{eq:DS1} and \eqref{eq:Ziju} that}
\begin{eqnarray}
  \lefteqn{ \pk{\exists_{(s,t)\in\Delta_{i,u}\times \Delta_{j,u}} (\vk{X}(s)+\vk{X}(t))_I>\sqrt{u}(2\vk{\alpha}_I+(s+t)\vk{\mu}_I)}
  }\nonumber\\
&=&\pk{\exists_{(s,t)\in[0,T]^2}
    ( \vk Z_{i,j;u})_I+\frac{1}{2\sqrt{u}}  (\vk{X}_1(s)+\vk{X}_1(T)+
    \vk{X}_2(t)-(s+t)\vk{\mu})_I  >
      \sqrt{u}(\vk{b}_{i,j;u})_I}\nonumber\\
    &=&
\frac{u^{-{m}/{2}}}{\sqrt{(2\pi)^{m}|(\Sigma_{i,j;u})_{II}|}}
\int_{\R^{m}}
\exp\LT(-\frac{1}{2}(\sqrt{u}(\vk{b}_{i,j;u})_I-\frac{\vk{x}_I}{\sqrt{u}})^\top
(\Sigma_{{i,j;u}})_{II}^{-1}(\sqrt{u}(\vk{b}_{i,j;u})_I-\frac{\vk{x}_I}{\sqrt{u}})\RT)\nonumber\\
&&\times
\pk{\exists_{(s,t)\in[0,T]^2}
    \frac{1}{2}  (\vk{X}_1(s)+\vk{X}_1(T)+
    \vk{X}_2(t)-(s+t)\vk{\mu})_I  >
      \vk x_I}
\,d\vk{x}_I.\label{eq:PL10} 
\end{eqnarray}
\COM{

\begin{eqnarray}
  \lefteqn{\pk{\exists_{s\in\Delta_{i,u}} \vk{X}(s)_I>\sqrt{u}(\vk{\alpha}_I+s\vk{\mu}_I),
    \exists_{t\in\Delta_{i,u}} \vk{X}(t)_I>\sqrt{u}(\vk{\alpha}+t\vk{\mu})}}\nonumber\\
&&\le
  \pk{\exists_{(s,t)\in\Delta_{i,u}\times \Delta_{j,u}} \vk{X}(s)_I+\vk{X}(t)_I>\sqrt{u}(2\vk{\alpha}_I+(s+t)\vk{\mu}_I)}\nonumber\\
  &&=
    \pk{\exists_{[0,\frac{T}{u}]\times [0,\frac{T}{u}]}
      \vk{X}(t_0+\frac{iT}{u}+s)+\vk{X}(t_0+\frac{jT}{u}+t)>
      \sqrt{u}(2\vk{\alpha}_I+(2t_0+\frac{i+j}{u}+s+t)\vk{\mu}_I)}\nonumber  \\
    &&= \pk{\exists_{[0,T]\times [0,T]}
     Z_{i,j;u}+\frac{(\vk{X}_1(s)+\vk{X}_1(T)+\vk{X}_2(t))_{I}-(s+t)\vk{\mu}}{2\sqrt{u}}    >
      \sqrt{u}(\vk{\alpha}_I+(t_0+\frac{(i+j)T}{2u})\vk{\mu}_I)}\nonumber\\
    &=&
u^{-\frac{m}{2}}
\int_{\R^{m}}
\frac{1}{(\sqrt{2\pi)^{m}|(\Sigma_{i,j;u})_{II}|}}
\exp(-\frac{1}{2}(\sqrt{u}(\vk{b}_{i,j;u})_I-\frac{\vk{x}_I}{\sqrt{u}})^\top
(\Sigma_{{i,j;u}})_{II}^{-1}(\sqrt{u}(\vk{b}_{i,j;u})_I-\frac{\vk{x}_I}{\sqrt{u}}))\nonumber\\
&&\times
\pk{ \exists_{s,t\in[0,T]}
    (\vk{X}_1(s)+\vk{X}_1(T)+\vk{X}_2(t))_I-(s+t)\vk{\mu})_I
   >\vk{x}_I,}\,d\vk{x}_I.\label{eq:PL10}
\end{eqnarray}
}
In particular for $i=0, j=2$
using \eqref{eq:PL10} and similar arguments as in the proof of Lemma \ref{lem:RPL}  we have
\begin{eqnarray}
 \lefteqn{\pk{\exists_{(s,t) \in \Delta_{0,u}\times \Delta_{2,u}}\ (\vk{X}(s)+\vk{X}(t))_I > \squ (2\vk{\alpha}_I+ \vk{\mu}_I (s+t) )}}\nonumber\\
  &&\sim  \widetilde{\HAS}(T)
\frac{  u^{-m/2} }{\sqrt{(2\pi t_0)^{m}\abs{\Sigma_{II}}}}
\exp(-\frac{u}{2}g_I(t_0))
\exp(-3aT),\label{eq:BB}
\end{eqnarray}
where
$a=\frac{g_I(t_0)}{8t_0}$ and
$$
\widetilde{\HAS}(T)= \int_{\R^{m}}e^{\frac{\x_I^\top \SI_{II}^{-1}\vk{b}_I}{ t_0}}
\pk{\exists_ {(s,t)\in [0, T ]^2}     \frac{1}{2}(
 \vk{X}_1(s)+\vk{X}_1(T)+\vk{X}_2(t) -  \vk{\mu}
 (s+ t))_I> \vk{x}_I}\, d\vk{x}_I.
$$
Note that in \EHb{a}  similar vein as in Lemma \ref{lem:HTT} 
we can prove the finiteness of $\widetilde{\HAS}(T)$.

We shall need an upper bound \EHb{for} $p_{i,j;u}$ \EHc{derived} in the following lemma.

\BEL\label{lemBB} For any fixed $T>0$, there exists some small $\vn>0$ such that, for all
$i,j$ satisfying $-N_u\le i <j\le N_u$,
\BQN\label{eq:upperpij}
p_{i,j;u}&\le&  CT^2 e^{2a  T}  u^{-m/2}
\exp\LT(-\frac{g_I(t_0)}{2} u\RT)\nonumber\\
&&\times\exp\LT(-  \frac{ g_I^{''}(t_0)-\vn}{4}\LT(\frac{ i T}{ \sqrt u} \RT)^2 \RT)\exp\Bigl(- a_\vn ((j-i+1)T)\Bigr)
\EQN
holds for some constant $C>0$ independent of $i,j,u$ and $T$, when $u$ is large,
where
$$
a_\vn = \frac{1}{ 2}  \LT(\frac{ g_I(t_0)-\vn }{4 t_0 +\vn} - \frac{\vn(g_I^{''}(t_0)-\vn)}{2} \RT),\ \ \ a=a_0=\frac{g_I(t_0)}{8t_0}.
$$
\EEL

{\bf Proof:}
. 
 In view of \eqref{eq:DS1} and \eqref{eq:PL10}  we have (\EHc{recall} $\vk{b}_{i,j;u}=\vk{b}(t_0+\frac{(i+j)T}{2u})$)
\begin{eqnarray}
p_{i,j;u} 
    &\le&
\frac{u^{-{m}/{2}}}{\sqrt{(2\pi)^{m}|(\Sigma_{i,j;u})_{II}|}}
\int_{\R^{m}}
\exp\LT(-\frac{1}{2}(\sqrt{u}(\vk{b}_{i,j;u})_I-\frac{\vk{x}_I}{\sqrt{u}})^\top
(\Sigma_{{i,j;u}})_{II}^{-1}(\sqrt{u}(\vk{b}_{i,j;u})_I-\frac{\vk{x}_I}{\sqrt{u}})\RT)\nonumber\\
&&\times
\pk{\exists_{(s,t)\in[0,T]^2}
    \frac{1}{2}  (\vk{X}_1(s)+\vk{X}_1(T)+
    \vk{X}_2(t)-\vk{\mu}(s+t))_I  >
      \vk x_I}
\,d\vk{x}_I.\label{eq:pijule}
\end{eqnarray}
Let for $T$ positive
$$
\widetilde{\H}_{I,i,j;u}(T)= \int_{\R^{m}}e^{\frac{\x_I^\top \SI_{II}^{-1}\vk{b}(t_0+\frac{(i+j)T}{2u})_I}{ t_0+\frac{j+3i-1}{4u}    }}
\pk{\exists _{(s,t)\in [0, T ]^2}     \frac{1}{2}(
\vk{X}_1(s)+\vk{X}_1(T)+\vk{X}_2(t)-  \vk{\mu}
 (s+ t))_I> \vk{x}_I}\, d\vk{x}_I.
$$
Since $\widetilde{\H}_{I,i,j;u}(T)\to \widetilde{\HAS}(T)$ as $u\to\infty$ uniformly with respect \EH{to} $-N_u\le i<j\le N_u$
we have that for large $u$
$$\widetilde{\H}_{I,i,j;u}(T)\le {\rm const}\ \widetilde{\HAS}(T).$$
\cc{Now for the expression  \EH{in} the exponent in
(\ref{eq:pijule}) we have
that
\BQNY
&&\Bigl(\sqrt{u}(\vk{b}_{i,j;u} )_I-\frac{\vk{x}_I}{\sqrt{u}}\Bigr)^\top
(\Sigma_{{i,j;u}})_{II}^{-1}\Bigl(\sqrt{u}(\vk{b}_{i,j;u})_I-\frac{\vk{x}_I}{\sqrt{u}}\Bigr)\\
&&=u \frac{1}{c_{i,j;u}} (\vk{b}_{i,j;u})_I^\top\Sigma_{II}^{-1} (\vk{b}_{i,j;u})_I - 2 \frac{1}{c_{i,j;u}} \vk x_I^\top\Sigma_{II}^{-1} (\vk{b}_{i,j;u})_I  + \frac{1}{u c_{i,j;u}} \vk{x}_I^\top\Sigma_{II}^{-1} \vk{x}_I \\
 &&=  u \frac{t_0+\frac{(i+j)T}{2u}}{c_{i,j;u}}  g_I(t_0+\frac{(i+j)T}{\kk{2u}})  - 2 \frac{1}{c_{i,j;u}} \vk x_I^\top\Sigma_{II}^{-1} (\vk{b}_{i,j;u})_I  + \frac{1}{u c_{i,j;u}} \vk{x}_I^\top\Sigma_{II}^{-1} \vk{x}_I.
\EQNY }
 It follows that  for all $u>0$ large
\BQN
p_{i,j;u}\le  C_0  \widetilde{\HAS}(T)  \frac{  u^{-m/2} }{\sqrt{(2\pi (t_0-\vn))^{m}\abs{\Sigma_{II}}}}
\exp\LT(-\frac{u}{2} \LT(1+\frac{(j-i+1)T}{(4t_0+\frac{(3i+j-1)T}{u})u}\RT)
g_I\Bigl(t_0+\frac{(i+j)T}{2u}\Bigr)\RT)  
\EQN
for some $C_0>0$ and some small $\vn>0$.
Furthermore, we have that for the small $\vn$
 \BQNY
\frac{(3i+j-1)T}{u} < \vn
 \EQNY
and
\BQNY
g_I(t_0+\frac{(i+j)T}{2u})\ge g_I(t_0)+\frac12 (g_I^{''}(t_0)-\vn) \LT(\frac{(i+j)T}{2u}\RT)^2
\EQNY
for all $-N(u)\le i<j\le N(u)$ and $u$ large. Moreover, for any $j>i$
\BQNY
\LT(\frac{(i+j)T}{2u}\RT)^2= \LT(\frac{(j-i)T}{2u} +\frac{ i T}{ u} \RT)^2 \ge \LT(\frac{ iT}{ u} \RT)^2 + \frac{(j-i)T (iT)}{ u^2} \ge \LT(\frac{ iT}{ u} \RT)^2 - \vn\frac{(j-i+1)T  }{ u }
\EQNY
holds for all $u$ large. \EHb{Consequently}, for any $j>i$
\begin{eqnarray*}
\lefteqn{
\exp\LT(-\frac{u}{2} \LT(1+\frac{(j-i+1)T}{(4t_0+\frac{(3i+j-1)T}{u})u}\RT)
g(t_0+\frac{(i+j)T}{2u})\RT)}\\
&&\le \exp\LT(-\frac{u}{2} \LT(1+\frac{(j-i+1)T}{(4 t_0 +\vn) u}\RT) \LT(g_I(t_0)+ \frac12 (g_I^{''}(t_0)-\vn)\LT(\LT(\frac{ i T}{ u} \RT)^2 - \vn\frac{(j-i+1)T  }{ u }\RT)  \RT)\RT)\\
&&=\exp\LT(-\frac{g_I(t_0)}{2} u-  \frac{ g_I^{''}(t_0)-\vn}{4}\LT(\frac{ i T}{ \sqrt u} \RT)^2
-\frac{(j-i+1)T}{ 2}  \LT(\frac{ g_I(t_0)+ \frac{g_I^{''}(t_0)-\vn}{2} \LT(\LT(\frac{ i T}{ u} \RT)^2 - \vn\frac{(j-i+1)T  }{ u }\RT)   }{4 t_0 +\vn} - \frac{\vn(g_I^{''}(t_0)-\vn)}{2} \RT)
\RT).
\end{eqnarray*}
With the small given \EHb{positive}  $\vn$, for all $-N(u)\le i<j\le N(u)$ and all large $u$ we have
$$
 g_I(t_0)+ \frac{g_I^{''}(t_0)-\vn}{2} \LT(\LT(\frac{ i T}{ u} \RT)^2 - \vn\frac{(j-i+1)T  }{ u }\RT)  \ge  g_I(t_0)-\vn.
$$
Consequently,
\BQNY 
\lefteqn{
\exp\LT(-\frac{u}{2} \LT(1+\frac{(j-i+1)T}{(4t_0+\frac{(3i+j-1)T}{u})u}\RT)
g\Bigl(t_0+\frac{(i+j)T}{2u} \Bigr)\RT)}\nonumber\\
&&\le  \exp\LT(-\frac{g_I(t_0)}{2} u\RT)\exp\LT(-  \frac{ g_I^{''}(t_0)-\vn}{4}\LT(\frac{ i T}{ \sqrt u} \RT)^2 \RT)\exp\Bigl(- a_\vn ((j-i+1)T)\Bigr)
\EQNY
from which we \EH{obtain that}
\begin{eqnarray*}
p_{i,j;u}&\le&  C\  \widetilde{\HAS}(T)\frac{  u^{-m/2} }{\sqrt{(2\pi (t_0-\vn))^{m}\abs{\Sigma_{II}}}}
\exp\LT(-\frac{g_I(t_0)}{2} u\RT)\nonumber\\
&&\times\exp\LT(-  \frac{ g_I^{''}(t_0)-\vn}{4}\LT(\frac{ i T}{ \sqrt u} \RT)^2 \RT)\exp\Bigl(- a_\vn ((j-i+1)T)\Bigr)
\end{eqnarray*}
holds when $u$ is large.
\cc{Next, in order to complete the proof it is sufficient to show that for any positive integer $T$}
\BQNY
\wHAS(T)\le T^2 e^{2a T} \wHAS(1).
\EQNY
\EHb{For $T,u$ positive define}
\BQNY
&&E_u(T)=[t_0,t_0+Tu^{-1}]\times  [t_0+2Tu^{-1},t_0+3Tu^{-1}],\\
&&E_u(k,l)=[t_0+ku^{-1},t_0+(k+1)u^{-1}]\times  [t_0+(2T+l)u^{-1},t_0+(2T+l+1)u^{-1}],\\ &&\qquad\qquad \qquad   k,l=0,1, \cdots, T-1.
\EQNY
It follow from \eqref{eq:BB} that
\BQNY
\lim_{u\to\IF}\frac{ \pk{\exists _{(s,t) \in E_u(T)}  (\vk{X}(s)+  \vk{X}(t))_I  >   \squ (2\vk{\alpha}_I+ \vk{\mu}_I (s+t) )}}{ \frac{  u^{-m/2} }{\sqrt{(2\pi t_0)^{m}\abs{\Sigma_{II}}}}  \exp\LT(-\frac{u}{2}   g_I(t_0)\RT)}
=e^{-3aT}   \widetilde{\HAS}(T).
\EQNY
Similarly, we can show that 
\BQNY
 \lim_{u\to\IF}\frac{ \pk{\exists_ {(t,w) \in E_u(k,l)}   (\vk{X}(s)+  \vk{X}(t) )_I >   \squ (2\vk{\alpha}_I+ \vk{\mu} _I(s+t) )}}{ \frac{  u^{-m/2} }{\sqrt{(2\pi t_0)^{m}\abs{\Sigma_{II}}}}  \exp\LT(-\frac{u}{2}   g_I(t_0)\RT)}
 =e^{- a(2T+l-k+1)}   \widetilde{\HAS}(1).
\EQNY
Furthermore, since
$$E_u(T)\subset \cup_{k=0}^{T-1}\cup_{l=0}^{T-1} E_u(k,l)$$
 we obtain from the above two equalities that
\BQNY
e^{-3aT}   \widetilde{\HAS}(T)\le \sum_{k=0}^{T-1}\sum_{l=0}^{T-1} e^{- a(2T+l-k+1)}   \widetilde{\HAS}(1),
\EQNY
which yields that
\BQNY
   \widetilde{\HAS}(T)&\le& e^{aT}\sum_{k=0}^{T-1}\sum_{l=0}^{T-1} e^{- a(l-k+1)}   \widetilde{\HAS}(1)\\
   &&\le e^{2aT}T^2 \widetilde{\HAS}(1)
\EQNY
establishing the proof.
\QED

Now, we are ready to show \eqref{eq:Pip}.  Note that
$$
\Pi(u)=\sum_{-N_u\le i<j\le N_u}p_{i,j;u}=\underset{j=i+1}{\sum_{-N_u\le i<j\le N_u}p_{i,j;u}}+
\underset{j>i+1}{\sum_{-N_u\le i<j\le N_u}p_{i,j;u}}=:\Pi_1(u)+\Pi_2(u).
$$
For $\Pi_1(u)$ we have
\BQNY
\Pi_1(u)&=&\sum_{i=-N(u)}^{N(u)} \Bigg(\pk{\exists_{t\in\Del_{i;u}} \ \ \vk{X}(t)>\sqrt u (\vk{\alpha}+\vk{\mu} t) } +\pk{\exists_{t\in \Del_{(i+1);u} } \ \ \vk{X}(t)>\sqrt u (\vk{\alpha}+\vk{\mu} t) }\\
&&\ \ \ \ -
\pk{\exists_{t\in \Del_{i;u}\cup \Del_{(i+1);u}}  \ \ \vk{X}(t)>\sqrt u (\vk{\alpha}+\vk{\mu} t) }
\Bigg)=: S_1(u)+S_2(u)-S_3(u).
\EQNY
Recall that we have proved in \nelem{lem:HT} \ccj{and \nelem{Lem2}} that
\BQN \label{recoH} \lim_{T\to\IF}T^{-1}\HAS(T)=\H_I\in(0,\IF),
\EQN
hence using similar arguments as for \eqref{eq:PL0} to $S_i(u), i=1,2,3,$ we conclude that
\BQN\label{eq:doubsum1}
&&\lim_{T\to\IF} \lim_{u\to\IF}\frac{\Pi_1(u)}{u^{(1-m)/2} \exp\LT(-\frac{ g_I(t_0)}{2} u\RT) }\nonumber \\
&&\qquad=  \frac{1}{\sqrt{(2\pi t_0)^{m} \abs{\Sigma_{II}}}}   \int_{-\infty}^\infty e^{-\frac{g_I^{''}(t_0) x^2}{4}}\psi(x)\,dx   \lim_{T\to\IF}\LT(\frac{2\HAS(T)}{T}- \frac{\HAS(2 T)}{T}\RT)\nonumber\\
&&\qquad=0.
\EQN
For $\Pi_2(u)$ we have from \eqref{eq:upperpij} that, there exists some $\vn>0,$ such that
\begin{eqnarray*}
\Pi_2(u)
&\le &
 C  T e^{2aT} u^{\frac{1-m}{2}}
  \exp\Bigl(-\frac{u}{2}g_I(t_0)\Bigr) \\
&&\times \frac{T}{\sqrt{u}}  \sum_{-N_u\le i \le N_u}\exp\LT(-  \frac{ g_I^{''}(t_0)-\vn}{4}\LT(\frac{ i T}{ \sqrt u} \RT)^2 \RT) \sum_{j\ge 1} \exp\Bigl(- a_\vn (jT)\Bigr) \exp(- 2a_\vn T)
\end{eqnarray*}
holds for all large $u$ with some  $ C>0$, implying thus
\BQNY
\lim_{T\to\IF}\lim_{\vn\to 0}\lim_{u\to\IF}\frac{\Pi_2(u)}{u^{(1-m)/2} \exp\LT(-\frac{ g_I(t_0)}{2} u\RT) }=0,
\EQNY
\K1{which establishes} \eqref{eq:Pip}.

\prooftheo{Thm1}   First note that \EHb{the} finiteness of $\H_I$ is established in \nelem{lem:HT} and the lower bound is obtained in  \nelem{Lem2}.  
Furthermore, in view of \eqref{eq:thetaT}, \eqref{eq:PL0}, \eqref{eq:Pip} \Lc{and} 
\EHb{letting $T \to \IF$
we obtain \EHb{(recall \eqref{recoH})}
 $$
 p(u)\sim 
 \frac{1}{\sqrt{(2\pi t_0)^{m}|\Sigma_{II}|}  }  \H_I u^{\frac{1- m}{2}}e^{-\frac{u}{2}g_I(t_0)}\int_{\R} e^{-\frac{g_I^{''}(t_0) x^2}{4}}\psi(x)\,dx
 , \quad u\to\IF.
 $$
Moreover, by  \nelem{PitI}
$$ r(u)= o( p(u)), \quad u\to \IF.  $$
Consequently, the claim follows from \eqref{eq:pr}. } \QED

\def\cu{c_{\lambda,u}}

\prooftheo{KorrRT}
Define
\BQNY
\widehat{\tau}_u=\inf\{t\ge0: \vk X(t)>(\vk \alpha +\vk \mu t)\sqrt u\}.
\EQNY
Since $\tau_u=u\widehat{\tau}_u $, for any $s\in\R$
\begin{eqnarray*}
\pk{\frac{\tau_u -t_0 u}{\sqrt u} \le s  \big\lvert \tau_u<\IF}&=&
\frac
{\pk{\frac{\tau_u -t_0 u}{\sqrt u} \le s ,\tau_u<\IF}}
{\pk{\tau_u<\infty}}
\\
&=&\frac{\pk{u\widehat{\tau}_u\le ut_0+\sqrt{u}s}}{P(u)}\\
&=&\frac{\pk{\exists_{t\in[0, t_0+s/\sqrt u]} \vk X(t)>(\vk \alpha +\vk \mu t)\sqrt u}}{P(u)} .
\end{eqnarray*}
Using the same arguments as in the proof of Theorem \ref{Thm1}, we have
\BQNY
\pk{\exists_{t\in[0, t_0+s/\sqrt u]} \vk X(t)>(\vk \alpha +\vk \mu t)\sqrt u}&\sim&
\pk{\exists_{t\in[t_0-\ln(u)/\sqrt{u}, t_0+s/\sqrt u]} \vk X(t)>(\vk \alpha +\vk \mu t)\sqrt u}\\
&\sim&
\frac{\HAS}{\sqrt{(2\pi t_0)^m \abs{\Sigma_{II}}}} \int_{-\IF}^s e^{-\ggt \frac{x^2}{4}} \psi(x)\,dx \,  u^{\EH{\frac{1- m}{2}}}e^{- \frac{\gt}{2} u},\ \ u\to\IF.
\EQNY
\EHb{In order} to derive the above result the only \EHb{required} modification in the proof of \netheo{Thm1} is the replacement of $\sum_{-N_u-1\le j\le N_u}$  by $\sum_{-N_u-1\le j\le \lfloor \sqrt u s /T\rfloor}$ in $R_T(u)$, see \eqref{eq:rTu}.
Consequently, the claim follows and thus the proof is complete. \QED

\section{Appendix}

\subsection{Quadratic programming problem}\label{ss:QPP}
This subsection is concerned \K1{with} discussions on Lemma \ref{AL}, which will be useful for the analysis of the function $g$ in the next subsection.
Recall from Lemma \ref{AL}, \K1{that} $\widetilde{\vk b}$ is the optimal solution of   the quadratic programming problem $P_{M}(\vk b)$ with the \EH{essential} index set $I$.  Next, we define \cc{for $I^c\neq \emptyset$}
\BQN \label{eq:WW}
K=\{j\in I^c: \vk  b_j=\widetilde{\vk b}_j\}.
\EQN

We start with some important remarks on Lemma \ref{AL}. 
\begin{remark}\label{rem2dimB}
 i)  \EH{If there is a unique  index set $I$  with maximal number of elements such that \eqref{eq:IJi} holds,
then $I$ is the \EH{essential} index set of $P_{\Sigma}(\vk b)$. Otherwise, if there are  $I_1 \ldot I_l$ index sets
 which have the same maximal number of elements such that \eqref{eq:IJi} holds, then the unique \EH{essential} index set say $I= I_k$ satisfies additionally \eqref{eq:hii}.}

ii) Note that, for any $I_1$ satisfying $\{1,\ldots, d\}\supseteq I_1\supset  I$, $\widetilde{\vk b_{I_1}}=\widetilde{\vk b}_{I_1}$ is the unique solution of the quadratic programming problem $ P_{\Sigma_{I_1I_1}}(\vk b_{I_1})$. \cc{
If  further}  $$\Sigma_{I_1I_1}^{-1} \vk b_{I_1}\ge \vk0_{I_1}$$
holds, then $\widetilde{\vk b_{I_1}}=\vk  b_{I_1}$ and
$$
\Sigma_{I_1I_1}^{-1} \vk b_{I_1}=\LT(\begin{array}{ccc}
\Sigma_{II}^{-1} \vk b_I  \\
\vk 0_{I_1\setminus I}
\end{array}\RT),
$$
with
$\Sigma_{II}^{-1} \vk b_I>\vk0_I$; see also the proof of 1 of  Proposition 2.5 in \cite{ENJH02}.

iii)
Consider the case $d=2$ and let $\vk{b}$ with $ b_1=1,  b_2= b \in (-\infty,1]$.
Suppose for simplicity that $\Sigma$ is a correlation matrix
with $\sigma_{12}=\rho\in (-1,1)$. If $ b> \rho$, then $\widetilde {\vk{b}}= \vk{b}$ and thus
$I=\{1, 2\}$.
  If $ b=\rho$, then $I=\{1\}, K=\{2\}$.
Finally, for $ b< \rho$ we have $I=\{1\}, K=\emptyset$.
\end{remark}

\def\veta{\vk b}
\def\wveta{\widetilde {\vk{b}}}
\BEL \label{Lem:IL}
Let $I\subset \{1,\ldots, d\}$ be the \EH{essential} index set of the quadratic programming problem $P_{M}(\vk  b), \EH{\vk  b \in \R^d \setminus (-\IF,0]^d} $.
We have, for any $I_1$ satisfying $\{1,\ldots, d\} \supseteq I_1\supset  I$, if
\cc{ $\vk  b_{I_1}^\top M_{I_1I_1}^{-1}\vk b_{I_1}=\vk b_{I}^\top M^{-1}_{II}\vk b_{I}$} holds, then $I_1\subseteq I\cup K,$ with $K$ given by \eqref{eq:WW}.
\EEL
{\bf Proof:} \LL{Note that from Remark \ref{rem2dimB} ii) we have $\wveta_{I_1}=\widetilde{\vk b_{I_1}}$.}
In the light of \eqref{eq:new},
$$
\veta_{I_1}^\top M_{I_1I_1}^{-1}\wveta_{I_1}= \kk{\veta}_{I_1}^\top M_{I_1I_1}^{-1}\LL{\widetilde{\vk b_{I_1}}}=\veta_{I}^\top M_{II}^{-1}\veta_{I} $$
 and
 $$\wveta_{I_1}^\top M_{I_1I_1}^{-1}(\veta_{I_1}-\wveta_{I_1})= \LL{  (\veta_{I_1}-\wveta_{I_1})^\top M_{I_1I_1}^{-1} \widetilde{\vk b_{I_1}}} =(\veta_{I_1}-\wveta_{I_1})_{I}^\top M_{II}^{-1}\veta_{I} =0.
$$
Further, since
\BQNY
\veta_{\kk{I_1}}^\top \cc{M}_{I_1I_1}^{-1}\veta_{I_1}&=& \veta_{I_1}^\top M_{I_1I_1}^{-1}(\veta_{I_1}-\wveta_{I_1})+
\veta_{I_1}^\top M_{I_1I_1}^{-1}\wveta_{I_1}\\
&=& (\veta_{I_1}-\wveta_{I_1})^\top M_{I_1I_1}^{-1}(\veta_{I_1}-\wveta_{I_1})+ \wveta_{I_1}^\top M_{I_1I_1}^{-1}(\veta_{I_1}-\wveta_{I_1}) +\veta_{I}^\top M_{II}^{-1}\veta_{I}\\
&=& (\veta_{I_1}-\wveta_{I_1})^\top M_{I_1I_1}^{-1}(\veta_{I_1}-\wveta_{I_1}) +\veta_{I}^\top M_{II}^{-1}\veta_{I}
\EQNY
we obtain
$$
 (\veta_{I_1}-\wveta_{I_1})^\top M_{I_1I_1}^{-1}(\kk{\veta}_{I_1}-\wveta_{I_1})=0
$$
thus $
\veta_{I_1}=\wveta_{I_1}
$,
implying that $I_1\subseteq I\cup K$ and \EHb{hence} the proof is complete. \QED

 \subsection{Analysis of $g$} \label{sectggg}
 In this subsection we analyze \EHc{the} \K1{function}
$$g(t)=\frac{1}{t} \inf_{\vk{v} \ge \vk{\alpha}+\vk{\mu} t}  \vk{v}^\top \Sigma^{-1}  \vk{v}$$
\EH{defined already  in the Introduction}. In the sequel we will denote by $I(t)$
 the \EH{essential}  index set of the quadratic programming problem
$P_{\Sigma}(\vk{\alpha}+t\vk{\mu})$. If $I(t)^c\neq \emptyset$ \K1{we} define
\BQNY
K(t)&=&\{j\in I(t)^c:  \Sigma_{jI(t)}\Sigma_{I(t)I(t)}^{-1}(\vk\alpha +\vk \mu t)_{I(t)}=(\vk\alpha +\vk \mu t)_j\},\\
J(t)&=&\{j\in I(t)^c:  \Sigma_{jI(t)}\Sigma_{I(t)I(t)}^{-1}(\vk\alpha +\vk \mu t)_{I(t)}>(\vk\alpha +\vk \mu t)_j\}.
\EQNY
 Note that, when analysing \EH{the} function $g$, \EH{the} index set $K(t)$ plays the role of $K$ from Section \ref{ss:QPP}.

\COM{The following elementary result is needed for the  continuity of  $ g$, see \nelem{lem:gc}.
\BEL  \label{lem:gcont}
Let $\mathcal X,\mathcal Y$ be two topological spaces and $\mathcal Y$ is compact.
Suppose that $f:\mathcal X\times \mathcal Y\to\R$ is continuous. Then
$$f_1(x)=\inf_{y\in \mathcal Y}f(x,y)$$
is continuous on $\mathcal X$.
\EEL
}
\BEL\label{lem:gc}
We have  $g\in C(0,\IF)$.
\EEL

{\bf Proof:} Let $h(t)=g(t) t$. For $g\in C(0,\IF)$ it is sufficient that  $h\in C(0,\IF)$.
In view of \nelem{AL} we have that for any $t\ge0$ there exists some $\vk{v}^*_t$,  so that
$$
h(t)=(\vk{v}^*_t+\vk{\alpha}+\vk{\mu} t)^\top \Sigma^{-1} (\vk{v}^*_t+\vk{\alpha}+\vk{\mu} t) , 
$$
where
$${\vk{v}^*_t} =
\left(
\begin{array}{ccc}
 \vk 0_{I(t)}  \\
\Sigma_{I(t)^c I(t)}\Sigma_{I(t)I(t)}^{-1}   (\vk\alpha+\vk \mu t) _{I(t)}  -  (\vk\alpha+\vk \mu t) _{\kk{I(t)^c}}
\end{array}
\right)
.$$
 For any fixed $t_1\in(0,\IF)$, it is easy to see that in a neighbourhood of $t_1$, say $(t_1-\vn,t_1+\vn)$, with some small $\vn>0$, we have
\BQN
h(t)= \inf_{\sup_{t\in (t_1-\vn,t_1+\vn)} \vk{v}^*_t \ge \vk{v} \ge \vk{0}}  (\vk{v}+\vk{\alpha}+\vk{\mu} t)^\top \Sigma^{-1} ( \vk{v}+\vk{\alpha}+\vk{\mu} t), \ \ \ t\in (t_1-\vn,t_1+\vn).
\EQN
\EH{Since for two \cc{topological spaces}  $\mathcal X,\mathcal Y$ with $\mathcal Y$  compact we have
$$f(x)=\inf_{y\in \mathcal Y} q(x,y)$$
is continuous on $\mathcal X$, provided that $q:\mathcal X\times \mathcal Y\to\R$ is continuous,
 we immediately get that}
%
$h\in C(t_1-\vn,t_1+\vn)$. Consequently, $h\in C(0,\IF)$ follows since $t_1$ was chosen arbitrarily. \QED

We show next  that
$$I(t)=\sum_{j} I_j\mathbb{I}(t\in U_j),$$
where $\mathbb{I}(\cdot)$ is the indicator function and $U_j$'s are of the following form
 \begin{equation}\label{eq:intervals}
 (a_k,b_k), [a_k,b_k), (a_k,b_k], [a_k,b_k], \{a_k\}, (b_k,\IF), [b_k,\IF),
 \end{equation}
where $0<a_k<b_k< \infty$
and $I_j\subseteq \{1,\ldots,d\}$.
Since point intervals are theoretically possible, we call such \K1{function} {\it almost piecewise constant set} function.

 \BEL\label{lem:II}
$I(t), t\ge 0$ is an  almost piecewise constant set function.
\EEL
{\bf Proof:} First, by \nelem{AL} for any $t\ge0$ there exists a unique $I(t)$ satisfying
\BQN
&&\Sigma_{I(t)I(t)}^{-1} (\vk\alpha+\vk\mu t)_{I(t)}>\vk 0_{I(t)}, \label{inA}\\
&&\Sigma_{I(t)^c I(t)}\Sigma_{I(t)I(t)}^{-1} (\vk\alpha+\vk\mu t)_{I(t)} \ge (\vk\alpha+\vk\mu t)_{I(t)^c },\ \ \mathrm{if \ } I(t)^c\neq \emptyset. \label{inB}
\EQN
\EH{Next, for each  \cc{$V_k\subseteq \{1, \ldots,d\}$} we solve \eqref{inA} and \eqref{inB} with $I(t)$ substituted by $V_k$ and $I(t)^c$ substituted by $V^c_k=\{1, \ldots,d\}\setminus V_k$}.
\COM{ we solve the following inequalities for $t\ge 0$
\BQNY
&&\Sigma_{V_kV_k}^{-1} (\vk\alpha+\vk\mu t)_{V_k}>\vk 0_{V_k}, \quad
\Sigma_{V^c_kV_k}\Sigma_{V_kV_k}^{-1} (\vk\alpha+\vk\mu t)_{V_k} \ge (\vk\alpha+\vk\mu t)_{V^c_k},\ \ \mathrm{if \ } V_k^c\neq \emptyset.
\EQNY
}
 Since for each $V_k$ the solution is a convex set, by the linearity the solution \cc{(if it exists)} is in one of the following forms
 \BQNY
 (a_k,b_k), [a_k,b_k), (a_k,b_k], [a_k,b_k], \{a_k\}, (b_k,\IF), [b_k,\IF).
 \EQNY
Therefore, there exists some finite partition $\{U_1, \ldots, U_q\}$ of $[0,\IF)$, with
$$q\le \sum_{i=1}^d
{\small \small\small\small\small \left(
\begin{array}{ccc}
d   \\
i
\end{array}
\right)}
$$ some constant and
 $U_j$  an interval such that the index set $I(t)=I_j\subseteq \{1,\ldots,d\}$ for all $t\in U_j^o$, hence the proof is complete. \QED

\BEL \label{lem:tj}
For the boundary points $t_j=\overline{U_j}\cap \overline{U_{j+1}}, j=1\ldots,q-1$, \cc{we have $I(t_j)\subset \{1,\ldots,d\}$ and $K(t_j)\neq \emptyset$. Moreover,} $I(t_j)\subseteq I(t) \subseteq  I(t_j)\cup K(t_j)$ for all $t\in {U_j}\cup {U_{j+1}}$.
\EEL

{\bf Proof:} 
It follows from  Lemma \ref{AL} that \EH{\eqref{inA} holds for $t=t_j$.}
\COM{
	\BQNY
 \Sigma^{-1}_{I(t_j),I(t_j)} (\vk \alpha+\vk \mu t_j) _{I(t_j)}>\vk 0_{I(t_j)}.
\EQNY
}
By continuity, there exists some small $\delta>0$ such that for all $t\in(t_j-\delta,t_j+\delta)$
\BQNY
 \Sigma^{-1}_{I(t_j),I(t_j)} (\vk \alpha+\vk \mu t) _{I(t_j)}>\vk 0_{I(t_j)}.
\EQNY
  This implies that \EH{$I(t_j)$ has less than $d$ elements}, 
since otherwise we \K1{would} have $I(t)=\{1,\ldots,d\}$ for all $t\in \overline{U_j}\cap \overline{U_{j+1}}$, a contradiction with the fact that $t_j$ is a boundary point. Similarly, if  \eqref{inA} and \eqref{inB} holds for $t=t_j$ \cc{with $K(t_j)= \emptyset$},
\COM{\BQN
&& \Sigma^{-1}_{I(t_j),I(t_j)} (\vk \alpha+\vk \mu t_j) _{I(t_j)}>\vk 0_{I(t_j)},\\
&& \Sigma_{I(t_j)^c,I(t_j)} \Sigma^{-1}_{I(t_j),I(t_j)} (\vk \alpha+\vk \mu t_j) _{I(t_j)}>(\vk \alpha+\vk \mu t_j) _{I(t_j)^c},
\EQN
}
then by continuity we  conclude that $I(t)=I(t_j)$ for all $t\in \overline{U_j}\cap \overline{U_{j+1}}$, again a contradiction. Thus,  $K(t_j)\neq \emptyset$.

Now,  let $I(t)=I_{j+1}, t\in U_{j+1}^o$ and $I(t)=I_j, t\in U_{j}^o$. Without loss of generality, we only show $I(t_j)\subseteq I_{j}$ since $I(t_j)\subseteq I_{j+1}$ follows with the same arguments.  Notice that
\BQN \label{eq:Ktj}
&& \Sigma^{-1}_{I(t_j),I(t_j)} (\vk \alpha+\vk \mu t_j) _{I(t_j)}>\vk 0_{I(t_j)},\nonumber\\
&& \Sigma_{K(t_j),I(t_j)} \Sigma^{-1}_{I(t_j),I(t_j)} (\vk \alpha+\vk \mu t_j) _{I(t_j)}=(\vk \alpha+\vk \mu t_j) _{K(t_j)},\\
&& \Sigma_{J(t_j),I(t_j)} \Sigma^{-1}_{I(t_j),I(t_j)} (\vk \alpha+\vk \mu t_j) _{I(t_j)}>(\vk \alpha+\vk \mu t_j) _{J(t_j)}.\nonumber
\EQN
Since equations in \eqref{eq:Ktj} are linear \cc{in $t_j$ for fixed $I(t_j), K(t_j)$,} two cases will be distinguished.

\underline{Case 1}. It holds that
$$\Sigma_{K(t_j),I(t_j)} \Sigma^{-1}_{I(t_j),I(t_j)} (\vk \alpha+\vk \mu t) _{I(t_j)}>(\vk \alpha+\vk \mu t) _{K(t_j)}$$
  for all $t\in U_{j}^o$.

\underline{Case 2}. There exists  some index $i\in K(t_j)$   such that
$$\Sigma_{i,I(t_j)} \Sigma^{-1}_{I(t_j),I(t_j)} (\vk \alpha+\vk \mu t) _{I(t_j)}<(\vk \alpha+\vk \mu t) _{i}
$$
holds  for all $t\in U_{j}^o$.

For Case 1, by continuity we conclude that $I_j=I(t_j)$.
Next, we focus on Case 2, and show for this case $  I(t_j)\subset I_j$.

\def \hI {\hat I}
Denote $\hI=I(t_j)\cup \{i\}$. We can show that
\BQN \label{eq:hI}
 \Sigma^{-1}_{\hI,\hI} (\vk \alpha+\vk \mu t ) _{\hI}>\vk 0_{\hI}
\EQN
holds for all $t\in U_{j}^o$ such that $t-t_j$ is small, which, by Remark \ref{rem2dimB} i),  implies that
\BQN\label{eq:Ihat}
\sharp I_j \ge \sharp \hI = \sharp I(t_j)+1.
\EQN
 In fact, denoting $B= \Sigma^{-1}_{\hI,\hI}$ we have
\BQNY
 \Sigma^{-1}_{\hI,\hI} (\vk \alpha+\vk \mu t ) _{\hI}=
 \LT(\begin{array}{ccc}
B_{I(t_j)I(t_j)}(\vk \alpha+\vk \mu t ) _{I(t_j)} +B_{I(t_j),i}  (\vk \alpha+\vk \mu t ) _{i}\\
B_{i,I(t_j)}(\vk \alpha+\vk \mu t ) _{I(t_j)} +B_{i,i}  (\vk \alpha+\vk \mu t ) _{i}
\end{array}\RT).
\EQNY
Since $B$ is positive definite, $B_{i,i}>0$. By  \EHc{the} properties of block  positive definite matrix $B$, we have that
\BQNY
B_{I(t_j)I(t_j)}(\vk \alpha+\vk \mu t ) _{I(t_j)} +B_{I(t_j),i}  (\vk \alpha+\vk \mu t ) _{i}&=&\Sigma^{-1}_{I(t_j),I(t_j)} (\vk \alpha+\vk \mu t) _{I(t_j)} \\
&&+B_{I(t_j),i} \LT((\vk \alpha+\vk \mu t ) _{i}-\Sigma_{i,I(t_j)}\Sigma^{-1}_{I(t_j),I(t_j)} (\vk \alpha+\vk \mu t ) _{I(t_j)}  \RT)
\EQNY
and
\BQNY
B_{i,I(t_j)}(\vk \alpha+\vk \mu t ) _{I(t_j)} +B_{i,i}  (\vk \alpha+\vk \mu t ) _{i}= B_{i,i} \LT((\vk \alpha+\vk \mu t ) _{i}-\Sigma_{i,I(t_j)}\Sigma^{-1}_{I(t_j),I(t_j)} (\vk \alpha+\vk \mu t ) _{I(t_j)}  \RT)>0.
\EQNY
Then, since
$$\Sigma^{-1}_{I(t_j),I(t_j)} (\vk \alpha+\vk \mu t_j) _{I(t_j)}>\vk 0_{I(t_j)}, \quad (\vk \alpha+\vk \mu t_j ) _{i}=\Sigma_{i,I(t_j)}\Sigma^{-1}_{I(t_j),I(t_j)} (\vk \alpha+\vk \mu t_j ) _{I(t_j)} $$
 we conclude that \eqref{eq:hI} holds  for all $t\in U_{j}^o$ such that $t-t_j$ is small.

On the other hand, since $I(t)=I_j, t\in U^o_j$ we have
\BQN 
&& \Sigma^{-1}_{I_jI_j} (\vk \alpha+\vk \mu t ) _{I_j}>\vk 0_{I_j},\nonumber\\
&& \Sigma_{I_j^cI_j} \Sigma^{-1}_{I_jI_j} (\vk \alpha+\vk \mu t) _{I_j}>(\vk \alpha+\vk \mu t) _{I_j^c}\label{eq:SIIc}
\EQN
hold for all $t\in U^o_j$.
The reason why we do not have equality in \eqref{eq:SIIc} is that if for some row equality holds with some $t_1\in U^o_j$, then $I(t)=I_j, t\in U^o_j$ will be invalid by  linearity of the equation. \EH{Consequetly, letting}  $t\to t_j$ in the above inequalities we obtain
\BQN
&& \Sigma^{-1}_{I_jI_j} (\vk \alpha+\vk \mu t_j ) _{I_j}\ge \vk 0_{I_j},\label{eq:SII2} \\
&& \Sigma_{I_j^cI_j} \Sigma^{-1}_{I_jI_j} (\vk \alpha+\vk \mu t_j) _{I_j}\ge (\vk \alpha+\vk \mu t_j) _{I_j^c}.\label{eq:SIIc2}
\EQN
Suppose that the first $l$ rows (the corresponding index set is denoted by $\hI_1$) of $\Sigma^{-1}_{I_jI_j} (\vk \alpha+\vk \mu t_j ) _{I_j}$ are positive and the last $\sharp I_j-l$ rows (the corresponding index set is denoted by $\hI_2$)  are equal to 0.
Since   $I(t_j)$ is the \EH{essential} index set of $P_{\Sigma}(\vk \alpha+\vk \mu t_j )$, in view of Remark \ref{rem2dimB} i) we have $  l \le \sharp I(t_j)$.
Next, as in Remark \ref{rem2dimB} ii) (see also the proof of Proposition 2.5 in \cite{ENJH02})  we have
\BQNY
\Sigma^{-1}_{I_jI_j} (\vk \alpha+\vk \mu t_j ) _{I_j}=\LT(\begin{array}{ccc}
\Sigma_{\hI_1\hI_1}^{-1}  (\vk \alpha+\vk \mu t_j ) _{\hI_1} \\
  \vk 0_{\hI_2}
\end{array}\RT)
\EQNY
and
\BQN \label{eq:SIh12}
\Sigma_{\hI_2\hI_1}\Sigma_{\hI_1\hI_1}^{-1}  (\vk \alpha+\vk \mu t_j ) _{\hI_1}= (\vk \alpha+\vk \mu t_j )  _{\hI_2},\ \ \ \
\Sigma_{\hI_1\hI_1}^{-1}  (\vk \alpha+\vk \mu t_j ) _{\hI_1}>\vk 0_{\hI_1}.
\EQN
Then rewriting \eqref{eq:SIIc2}   we have
\BQNY
 \Sigma_{I_j^c \hI_1} \Sigma^{-1}_{\hI_1\hI_1} (\vk \alpha+\vk \mu t_j) _{\hI_1}\ge (\vk \alpha+\vk \mu t_j) _{I_j^c},
\EQNY
which together with \eqref{eq:SIh12} yields that $\hI_1$ is also an \EH{essential} index set of the problem $P_{\Sigma}(\vk \alpha+\vk \mu t_j )$. Thus, \cc{by uniqueness,} $I(t_j)=\hI_1\subset I_j$. Consequently, $  I(t_j)  \subseteq I(t) $ for all $t\in {U_j}\cup {U_{j+1}}$. Finally, we show $ I(t) \subseteq  I(t_j)\cup K(t_j)$ for all $t\in {U_j}\cup {U_{j+1}}$. Since $g\in C(0,\IF)$ we have
\BQNY
 (\vk \alpha+\vk \mu t_j) _{I(t_j)}^\top  \Sigma^{-1}_{I(t_j),I(t_j)} (\vk \alpha+\vk \mu t_j) _{I(t_j)}
& =& (\vk \alpha+\vk \mu t_j) _{I_j}^\top  \Sigma^{-1}_{I_j,I_j} (\vk \alpha+\vk \mu t_j) _{I_j}\\
& = &(\vk \alpha+\vk \mu t_j) _{I_{j+1}}^\top  \Sigma^{-1}_{I_{j+1},I_{j+1}} (\vk \alpha+\vk \mu t_j) _{I_{j+1}}.
\EQNY
Consequently, we conclude  from \nelem{Lem:IL} that $I_j\subseteq I(t_j)\cup K(t_j)$ and  $I_{j+1}\subseteq I(t_j)\cup K(t_j)$, establishing the proof. \QED

\COM{
If the above lemma cannot be proved then we make the following  hypothesis:

 {\bf H:} In a \e1{dimension-reduction}
instants $t^{'}$, if $I(t^{'})=I$ and $K(t^{'})=K$, then
either $I\cup K= I(t^{'}-)$ and    $ I(t^{'}+)=I$
or  $I\cup K= I(t^{'}+)$ and    $ I(t^{'}-)=I$.

The following elementary result is needed for the proof of \nelem{lem:gc1}.
\BEL  \label{lem:gcont}
Let $\mathcal X,\mathcal Y$ be two topological spaces and $\mathcal Y$ is compact.
Suppose that $f:\mathcal X\times \mathcal Y\to\R$ is continuous. Then
$$f_1(x)=\inf_{y\in \mathcal Y}f(x,y)$$
is continuous on $\mathcal X$.
\EEL

  \prooflem{lem:gc1}
 Let $h(t)=g(t) t$. For $g\in C^1(0,\IF)$ it is sufficient that  $h\in C^1(0,\IF)$.

\underline{Step 1: $h\in C(0,\IF)$.}
In view of \nelem{AL} we have, for any $t\ge0$ there exists some $\vk{v}^*_t$,  so that
$$
h(t)=(\vk{v}^*_t+\vk{\alpha}+\vk{\mu} t)^\top \Sigma^{-1} (\vk{v}^*_t+\vk{\alpha}+\vk{\mu} t) 
$$
where ${\vk{v}^*_t} =
\left(
\begin{array}{ccc}
 \vk 0_{I(t)}  \\
\Sigma_{I(t)^c I(t)}\Sigma_{I(t)I(t)}^{-1}   (\vk\alpha+\vk \mu t) _{I(t)}  -  (\vk\alpha+\vk \mu t) _{I(t)}
\end{array}
\right)
$. For any fixed $t_1\in(0,\IF)$, it is easy to see that in a neighbourhood of $t_1$, say $(t_1-\vn,t_1+\vn)$, with some small $\vn>0$, we have
$$
h(t)= \inf_{\sup_{t\in (t_1-\vn,t_1+\vn)} \vk{v}^*_t \ge \vk{v} \ge \vk{0}}  (\vk{v}+\vk{\alpha}+\vk{\mu} t)^\top \Sigma^{-1} ( \vk{v}+\vk{\alpha}+\vk{\mu} t), \ \ \ t\in (t_1-\vn,t_1+\vn).
$$
Thus, by \nelem{lem:gcont} we have $h\in C(t_1-\vn,t_1+\vn)$. Consequently, $h\in C(0,\IF)$ follows since $t_1$ was chosen arbitrarily.
}
\prooflem{lem:gc1}
By \nelem{lem:II} for any $ j=1 \ldot q$ we have
\BQNY
h(t)&=&  \inf_{\vk{v} \ge \vk{\alpha}+\vk{\mu} t}  \vk{v}^\top \Sigma^{-1}  \vk{v}\\
&=&   (\vk{\alpha}+\vk{\mu} t)^\top_{I(t)}  \Sigma^{-1}_{I(t),I(t)}  (\vk{\alpha}+\vk{\mu} t)_{I(t)}\\
&=&   \vk{\alpha} ^\top_{I(t)}  \Sigma^{-1}_{I(t),I(t)}   \vk{\alpha} _{I(t)} +
2 t \vk{\alpha} ^\top_{I(t)}  \Sigma^{-1}_{I(t),I(t)}   \vk{\mu}  _{I(t)}+
\vk{\mu}^\top_{I(t)}  \Sigma^{-1}_{I(t),I(t)}   \vk{\mu} _{I(t)} t ^2\\
&= & \vk{\alpha} ^\top_{I_j}  \Sigma^{-1}_{I_j,I_j}   \vk{\alpha} _{I_j} +
2 t \vk{\alpha} ^\top_{I_j}  \Sigma^{-1}_{I_j,I_j}   \vk{\mu}  _{I_j}+
\vk{\mu}^\top_{I_j}  \Sigma^{-1}_{I_j,I_j}   \vk{\mu} _{I_j} t^2, \ \ \ t\in U_j^o.
\EQNY
Clearly, $h\in C^1(U_j^o )$  for all  $ j=1 \ldot  q$. Thus, to prove that $g\in C^1(0,\IF)$ it is sufficient to show that, for any $t_j=\overline{U_j}\cap \overline{U_{j+1}}$,
$$\cc{h'(t_j+)=h'(t_j-)}$$
 holds.
 It follows that
 \BQNY
&& h'(t_j-)=2 ( \vk{\alpha} ^\top_{I_j}  \Sigma^{-1}_{I_j,I_j}   \vk{\mu}  _{I_j}+
\vk{\mu}^\top_{I_j}  \Sigma^{-1}_{I_j,I_j}   \vk{\mu} _{I_j} t_j)= 2 \vk{\mu}^\top_{I_j}  \Sigma^{-1}_{I_j,I_j}   (\vk \alpha+\vk{\mu}  t_j) _{I_j},\\
&& h'(t_j+) 
= 2 \vk{\mu}^\top_{I_{j+1}}  \Sigma^{-1}_{I_{j+1}I_{j+1}}   (\vk \alpha +\vk{\mu} t_j) _{I_{j+1}}.
 \EQNY
 Next, from \nelem{lem:tj}  we have $I(t_j)\subseteq I_j$ and $I(t_j)\subseteq I_{j+1}$.
For notational simplicity, we denote $B= \Sigma^{-1}_{I_j,I_j} $,  $J_j=I_j \setminus I(t_j)$. Since
$$
 \Sigma^{-1}_{I_j,I_j}   (\vk \alpha+\vk{\mu}  t) _{I_j}> \vk 0_{I_j},\ \ \ \ t\in U_j^o
$$
we have
$$
 \Sigma^{-1}_{I_j,I_j}   (\vk \alpha+\vk{\mu}  t_j) _{I_j}\ge \vk 0_{I_j}.
$$
Thus, by  Remark \ref{rem2dimB} ii)
\BQN\label{eq:Ijtj}
 \Sigma^{-1}_{I_j,I_j}   (\vk \alpha+\vk{\mu}  t_j) _{I_j}=
\left(
\begin{array}{ccc}
 \Sigma^{-1}_{I(t_j),I(t_j)}   (\vk \alpha+\vk{\mu}  t_j) _{I(t_j)} \\
\vk 0_{J_j}
\end{array}
\right),
\EQN
with $ \Sigma^{-1}_{I(t_j),I(t_j)}   (\vk \alpha+\vk{\mu}  t_j) _{I(t_j)}   >\vk 0_{I(t_j)}$  implying
$$
h'(t_j-)= 2 \vk{\mu}^\top_{I_j}  \Sigma^{-1}_{I_j,I_j}   (\vk \alpha+\vk{\mu}  t_j) _{I_j}=2  \vk{\mu}^\top_{I(t_j)}  \Sigma^{-1}_{I(t_j),I(t_j)}   (\vk \alpha+\vk{\mu}  t_j) _{I(t_j)}.
$$
Similarly, we have
$$
h'(t_j+)=2  \vk{\mu}^\top_{I(t_j)}  \Sigma^{-1}_{I(t_j),I(t_j)}   (\vk \alpha+\vk{\mu}  t_j) _{I(t_j)}.
$$
Consequently, $g\in C^1(0,\IF)$ is proved.

Now,
\BQNY
g(t)
=\frac{1}{t}    \vk{\alpha} ^\top_{I_j}  \Sigma^{-1}_{I_j,I_j}   \vk{\alpha} _{I_j} +
2 \vk{\alpha} ^\top_{I_j}  \Sigma^{-1}_{I_j,I_j}   \vk{\mu}  _{I_j}+
\vk{\mu}^\top_{I_j}  \Sigma^{-1}_{I_j,I_j}   \vk{\mu} _{I_j} t, \ \ \  t\in U_j^o
\EQNY
and
$$
g'(t)=\frac{\vk{\mu}^\top_{I_j}  \Sigma^{-1}_{I_j,I_j}   \vk{\mu} _{I_j} t^2 -  \vk{\alpha} ^\top_{I_j}  \Sigma^{-1}_{I_j,I_j}   \vk{\alpha} _{I_j} }{t^2}, \ \ \  t\in U_j^o.
$$
Since for any  nonempty  $I_j\subset\{1 \ldot d\}$
$$
\vk{\alpha} ^\top_{I_j}  \Sigma^{-1}_{I_j,I_j}   \vk{\alpha} _{I_j} >0,\ \ \vk{\mu}^\top_{I_j} \Sigma^{-1}_{I_j,I_j}   \vk{\mu} _{I_j}>0
$$
we have $g(t)\to\IF$ as $t\to \IF$ and $t\to0$, and $g'(t)<0$ for all $t$ around 0, $g'(t)>0$ for all $t$ large enough. Thus, \EH{the} function $g$ has \EH{a unique minimizer \cc{in} $[0,\IF]$}. Note that function $a/s +b+cs$ is decreasing to the left of some $s_0>0$ and increasing to the right.  Consider the interval $U_j$. The function $g$ has a unique minimum on  $U_j$. 
If at $t_j$ the function is decreasing it either decreasing
in the whole interval, or $t_0$ belongs to $U_j$ so it is increasing at
$t_{j+1}$ and consequently it is increasing at each entrance to constancy interval $U_k, k>j$.
In this case, \eqref{eq:t0} holds and  $g'(t_0)=0$. Next,  for $t_0$ we have from \nelem{lem:II} that there exist some small $\vn>0$ and $I^+,I^-\subseteq \{1 \ldot d\}$ such that
\BQNY
&&g(t)=\frac{1}{t}    \vk{\alpha} ^\top_{I^+}  \Sigma^{-1}_{I^+I^+}   \vk{\alpha} _{I^+} +
2 \vk{\alpha} ^\top_{I^+}  \Sigma^{-1}_{I^+I^+}   \vk{\mu}  _{I^+}+
\vk{\mu}^\top_{I^+}  \Sigma^{-1}_{I^+I^+}   \vk{\mu} _{I^+} t, \ \ \  t\in (t_0,t_0+\vn),\\
&&g(t)=\frac{1}{t}    \vk{\alpha} ^\top_{I^-}  \Sigma^{-1}_{I^-I^-}   \vk{\alpha} _{I^-} +
2 \vk{\alpha} ^\top_{I^-}  \Sigma^{-1}_{I^-I^-}   \vk{\mu}  _{I^-}+
\vk{\mu}^\top_{I^-}  \Sigma^{-1}_{I^-I^-}   \vk{\mu} _{I^-} t, \ \ \  t\in (t_0-\vn,t_0).
\EQNY
Then it follows that \eqref{eq:gt0pm} holds. \QED


\subsection{Analysis of 2-dimensional case}\label{ss.anal.2-dim}
We now demonstrate details for Section \ref{ss.two-dim}.
Recall that in our notation  $I(t)$ is the \EH{essential index} set of the quadradtic problem $P_\Sigma( \vk\alpha +\vk \mu t )$. If $I(t)^c\neq \emptyset$ \Lc{we} define
$$
K(t)=\{j\in I(t)^c:  \Sigma_{jI(t)}\Sigma_{I(t)I(t)}^{-1}(\vk\alpha +\vk \mu t)_{I(t)}=(\vk\alpha +\vk \mu t)_j\}.
$$
 Further define
$$
b_t=\frac{\alpha_2+t}{\alpha_1+t}\in (0,1).
$$
It follows that
\BQNY
g(t)=\frac{(\alpha_1+t)^2}{t}\inf_{\vk{v}\ge\vk{b}_t } \vk{v}^\top\Sigma^{-1}\vk{v},\ \ \ \vk{b}_t=(1,b_t)^\top.
 \EQNY
\underline{Case 1. $\rho<0$}. Clearly $b_t>\rho$ \EH{and thus in} view of \Lc{Remark \ref{rem2dimB} iii)}  we have  that $I(t)=\{1,2\},t>0$ and   
\BQNY
\inf_{\vk{v}\ge\vk{b}_t } \vk{v}^\top\Sigma^{-1}\vk{v} =\frac{1}{1-\rho^2}(1+b_t^2-2b_t \rho)
\EQNY
\EHb{implying}
\BQNY
g(t)=g_1(t):=\frac{(\alpha_1+t)^2 }{t} \frac{1}{1-\rho^2}(1+b_t^2-2b_t \rho).
\EQNY
Note that we slightly abuse the notation writing $g_1$ instead \cc{of $g_{\{1,2\}}$.} It \EH{follows} that for
$$t_0^{(1)}=\sqrt{\frac{\alpha_1^2+\alpha_2^2-2\alpha_1\alpha_2\rho}{2(1-\rho)}}\EHc{>0}$$
we have
$$\inf_{t\ge0}g(t)= g_1(t_0^{(1)})=\frac{2}{1+\rho}(\alpha_1+\alpha_2+2t_0^{(1)}).$$

 \underline{Case 2. $\rho>0$}. In such a case, we have to consider if $b_t>\rho$ or not. Several different sub-cases are thus discussed in the following.

 \underline{Case 2.1. $\alpha_1 \rho\le \alpha_2$}. For this case, we have always $b_t>\rho, t>0$. Then $I(t)=\{1,2\},t>0$ and $g(t)=g_1(t)$. 

 \underline{Case 2.2. $\alpha_1 \rho>\alpha_2$}. Let
 $$ Q:=\frac{\alpha_1\rho-\alpha_2}{1-\rho}.$$
 We have
 \BQNY
 &&(a)\ \ \  \{b_t>\rho\} \ \Leftrightarrow \ \{t> Q\},\ \  \mathrm{for\ which}\ \    I(t)=\{1,2\},    \\
 &&(b)\ \ \  \{b_t<\rho\} \ \Leftrightarrow \ \{t< Q\},\ \  \mathrm{for\ which}\ \    I(t)=\{1\}, \ K(t)=\emptyset,          \\
 &&(c)\ \ \  \{b_t=\rho\} \ \Leftrightarrow \ \{t= Q\}, \ \  \mathrm{for\ which}\ \    I(t)=\{1\}, \ K(t)=\{2\}.
 \EQNY
Now consider (a). Since $b_t>\rho$, we have $g(t)=g_1(t), t>Q$. Now we have to check if $\tzo>Q$ or not.
We can show that
\BQNY
\tzo>Q\ \  \Leftrightarrow \ \ \rho<\frac{\alpha_1+\alpha_2}{2\alpha_1}.
\EQNY
Thus, we have

(a1). If $\alpha_2/\alpha_1<\rho<\frac{\alpha_1+\alpha_2}{2\alpha_1}$, then $\inf_{t\in(Q,\IF)}g(t)=g_1(\tzo)$;

(a2). If $ \rho>\frac{\alpha_1+\alpha_2}{2\alpha_1}$ then $\inf_{t\in(Q,\IF)}g(t)=g_1(Q)$.

Next consider (b).
 Let $g_2(t)=(\alpha_1+t )^2/t$ which attains its minimum at the unique point $\tzt=\alpha_1$. Since  $b_t<\rho$, we have $g(t)=g_2(t), t\in[0,Q)$. Similarly as above we have to check if $\tzt<Q$. We can show that
 \BQNY
 \tzt<Q\ \ \Leftrightarrow\ \ \rho>\frac{\alpha_1+\alpha_2}{2\alpha_1}.
 \EQNY
Thus we have

(b1). If $\alpha_2/\alpha_1<\rho<\frac{\alpha_1+\alpha_2}{2\alpha_1}$, then $$\inf_{t\in[0,Q)}g(t)=g_2(Q)=\frac{(\alpha_1-\alpha_2)^2}{(1-\rho)(\alpha_1\rho-\alpha_2)}.$$

(b2). If $ \rho>\frac{\alpha_1+\alpha_2}{2\alpha_1}$, then
$$\inf_{t\in[0,Q)}g(t)=g_2(\tzt)=4\alpha_1.$$
Furthermore, by \EH{the} definitions of $g_1,g_2$ and $Q$ we obtain
\BQNY
g_1(Q)=g_2(Q).
\EQNY
The above findings are summarized in the following lemma:
\BEL \label{lem:two}
(1). If $ -1<\rho\le  \alpha_2/\alpha_1$, then
$I(t)= \{1,2\},t>0
$  and
\BQNY
t_0=\tzo,\ \ I= \{1,2\},\ \ g_I(t_0)=g_1(\tzo)  ,\ \  g_I^{''}(t_0)=g_1^{''}(\tzo)=2(\tzo)^{-3}\frac{\alpha_1^2+\alpha_2^2-2\alpha_1\alpha_2\rho}{1-\rho^2}.
\EQNY

(2). If $\alpha_2/\alpha_1<\rho<\frac{\alpha_1+\alpha_2}{2\alpha_1}$, then
\BQNY
I(t)= \{1\},\ \ 0<t\le Q, \ \ \ \ I(t)= \{1,2\},\ \ t>Q
\EQNY
and
\BQNY
t_0=\tzo>Q,\ \ I= \{1,2\},\ \ g_I (t_0) =g_1(\tzo) ,\ \   g_I^{''}(t_0)=g_1^{''}(\tzo).
\EQNY

(3). If $ \rho=\frac{\alpha_1+\alpha_2}{2\alpha_1}$, then
\BQNY
I(t)= \{1\},\ \ 0<t\le Q, \ \ \ \ I(t)= \{1,2\},\ \ t>Q
\EQNY
and
\BQNY
t_0=\tzo=\tzt=Q,\ \ I= \{1\},\ \ K=\{2\},\ \  g_I(t_0) =g_2(\tzt),\ \   g_I^{''}(t_0)=  g_2^{''}(\tzt)=2\alpha_1^{-1}.
\EQNY

(4).  If $ \frac{\alpha_1+\alpha_2}{2\alpha_1}<\rho<1$, then
\BQNY
I(t)= \{1\},\ \ 0<t\le Q, \ \ \ \ I(t)= \{1,2\},\ \ t>Q
\EQNY
and
\BQNY
t_0=\tzt<Q,\ \ I= \{1\},\ \ K=\emptyset,\ \  g_I(t_0) =g_2(\tzt) ,\ \   g_I^{''}(t_0)=2\alpha_1^{-1}.
\EQNY

\EEL

\begin{remark}  We point out that in general the second derivative of $g$ at $t_0$ is discontinuous. For instance, for the case where $\rho=\frac{\alpha_1+\alpha_2}{2\alpha_1}$ in \nelem{lem:two}  we have
\cc{
	\begin{eqnarray*}
		g(t)=\left\{
		\begin{array}{cc}
			\frac{1}{t} (\alpha_1+t)^2,& 0<t \le  \alpha_1,\\
			\frac{1}{t} (\vk \alpha+\vk 1 t)^\top \Sigma^{-1} (\vk \alpha+\vk 1 t), & t> \alpha_1.
		\end{array}
		\right.
	\end{eqnarray*}
	Hence}
\begin{eqnarray*}
	g^{'}(t)=\left\{
	\begin{array}{cc}
		1-\frac{\alpha_1^2}{t^2},&
		0<t \le  \alpha_1,\\
		\frac{2}{1+\rho}-\frac{4\alpha_1^3}{(3\alpha_1+\alpha_2)t^2},&t> \alpha_1
	\end{array}
	\right.
, \quad 	g^{''}(t)=\left\{
	\begin{array}{cc}
		\frac{2\alpha_1^2}{t^3},&
		0<t \le  \alpha_1,\\
		 \frac{8\alpha_1^3}{(3\alpha_1+\alpha_2)t^3},&t> \alpha_1.
	\end{array}
	\right.
\end{eqnarray*}
\e1{Consquetnly}, $g\in C^1(0,\IF)$ is  decreasing in the interval $(0,  \alpha_1)$. \e1{Its}
first derivative is 0 at  $t_0=\alpha_1$, however its second derivative is not continuous at $t_0$.

\end{remark}

\subsection{ Proof of \eqref{rT}} \label{TDT}
Recall $R_T(u)$ defined in \eqref{eq:bPL}.
We  \EH{derive next}  sharper bounds for $P_{j;u}(T,\vk x_I)$ and $ f_{j;u}(T,\vk x_I)$. Since  $\Sigma_{JI}\Sigma_{II}^{-1}\vk{b}_I<\vk{b}_J$, then for any small $\vn>0$ and any large $Q>0$
$$
P_{j;u}^-(T,\vk x_I,\vn,Q)\le P_{j;u}(T,\vk x_I)\le P_{j;u}^+(T,\vk x_I,\vn)
$$
holds for all $-N_u\le j\le N_u$ when $u$ is large enough, where
\BQNY
&&P_{j;u}^-(T,\vk x_I,\vn,Q)=\pk{
	\exists_{t\in[0,T]}
	\begin{array}{l}
		(\vk{X}(t)-t\vk{\mu})_I>\vk{x}_I\\
		\sqrt{t_0-\vn} Y_{K}
		-\vn\abs{\vk{Z}_K(t,\vk{x}_I) }>\frac{jT}{\sqrt{u}}
		(\vk{\mu}_K-\Sigma_{KI}\Sigma_{II}^{-1}\vk{\mu}_I)\\
		\sqrt{t_0-\vn} Y_{J}-\vn\abs{
			\vk{Z}_J(t,\vk{x}_I)}
		>-Q\vk 1_J
	\end{array}},  \\
	&&P_{j;u}^+(T,\vk x_I,\vn)=\pk{
		\exists_{t\in[0,T]}
		\begin{array}{l}
			(\vk{X}(t)-t\vk{\mu})_I>\vk{x}_I\\
			\sqrt{t_0+\vn} Y_{K}
			+\vn\abs{\vk{Z}_K(t,\vk{x}_I) }>\frac{jT}{\sqrt{u}}
			(\vk{\mu}_K-\Sigma_{KI}\Sigma_{II}^{-1}\vk{\mu}_I)
		\end{array}}.
		\EQNY
		Furthermore, for any large $L>0$, \EH{we can find $\vn>0$} sufficiently small such that
		\BQNY
		e^{  \frac{  \x _{I }^\top (\SI_{I I }^{-1}\b_I - \vk\vn_I^{\x_I})}{t_0-\vn(\x_I)} -\vn } \le
		f_u(T,j,\vk{x}_I) 
		\EQNY
		holds for all $\abs{\abs{\vk x_I}}\le L$ and all $-N_u\le j\le N_u$ when $u$ is large enough,
		where
		\BQN\label{eq:vnx1}
		\cc{\vn(\x_I)}=\left\{
		\begin{array}{ll}
			-\vn , &  \x _{I }^\top  \SI_{I I }^{-1}\b_I >0,\\
			\vn , &  \x _{I }^\top  \SI_{I I }^{-1}\b_I\le 0,
		\end{array}
		\right.  \EQN
		and \cc{$\vk{\vn}_I^{\x_I}=( \vn_i^{\x_I}, i\in I) $ with}
		\BQN\label{eq:vnx2}
		\cc{\vn_i^{\x_I}}=\left\{
		\begin{array}{ll}
			\vn , &  x_i >0,\\
			- \vn , & x_i\le 0,
		\end{array}
		\right.\ \ \ \ i\in I.
		\EQN
		Similarly,
		\BQNY
		f_u(T,j,\vk{x}_I)  \le e^{  \frac{  \x _{I }^\top (\SI_{I I }^{-1}\b_I + \vk\vn_I^{\x_I})}{t_0+\vn(\x_I)}  }
		\EQNY
		holds for all $\vk x_I\in \R^m$, $-N_u\le j\le N_u$ when $u$ is large enough.
		Moreover, it follows from \eqref{eq:gt0t1} that for the given $\vn$
		\BQNY
		\frac{g_I^{''}(t_0)-\vn}{2}\LT(\frac{jT}{  u}\RT)^2 \le g_I(t_0+\frac{jT}{  u})-g_I(t_0)\le \frac{g_I^{''}(t_0)+\vn}{2}\LT(\frac{jT}{  u}\RT)^2
		\EQNY
		holds for all $-N_u\le j\le N_u$ when $u$ is large enough.
		Consequently, we obtain the  following upper bound
		\BQNY
	\cc{	R_T(u)}\le  \frac{1}{(t_0-\vn)^{m/2}}  ( F_1(L,\vn,T, u)+ F_2(L,\vn,T, u)),
		\EQNY
		where
		\BQNY
		F_1(L,\vn,T, u)&=&\frac{T}{\sqrt{u}}\sum_{-N_u-1\le j\le N_u}
		\exp\LT(-\frac{g_I^{''}(t_0)-\vn }{4}  \LT(\frac{jT}{\sqrt u}\RT)^2     \RT)
		\int_{\abs{\abs{\vk x_I}}\le L}
		e^{  \frac{  \x _{I }^\top (\SI_{I I }^{-1}\b_I + \vk\vn_I^{\x_I})}{t_0+\vn(\x_I)}  } P_{j;u}^+(T,\vk x_I,\vn)\,d\vk{x}_I \\
		F_2(L,\vn,T, u)&=&\frac{T}{\sqrt{u}}\sum_{-N_u-1\le j\le N_u}
		\exp\LT(-\frac{g_I^{''}(t_0)-\vn }{4}  \LT(\frac{jT}{\sqrt u}\RT)^2     \RT)\\
		&&    \times   \int_{\abs{\abs{\vk x_I}}> L}
		e^{  \frac{  \x _{I }^\top (\SI_{I I }^{-1}\b_I + \vk\vn_I^{\x_I})}{t_0+\vn(\x_I)}  } \pk{\exists_{t\in[0,T]}(\vk{X}(t)-t\vk{\mu})_I>\vk{x}_I}\,d\vk{x}_I.
		\EQNY
		Next, it follows that
		\BQNY
		\lim_{\vn\to 0}\lim_{u\to0} F_1(L,\vn,T, u) &=&    \int_{\abs{\abs{\vk x_I}}\le L}
		e^{  \frac{  \x _{I }^\top \SI_{I I }^{-1}\b_I }{t_0}  }  \pk{\exists_{t\in[0,T]}(\vk{X}(t)-t\vk{\mu})_I>\vk{x}_I}\,d\vk{x}_I
		\int_{-\infty}^\infty e^{-\frac{g_I^{''}(t_0) x^2}{4}}\psi(x)\,dx
		\EQNY
		and
		\BQNY
		\lim_{\vn\to 0}\lim_{u\to0} \frac{T}{\sqrt{u}}\sum_{-N_u-1\le j\le N_u}
		\exp\LT(-\frac{g_I^{''}(t_0)-\vn }{4}  \LT(\frac{jT}{\sqrt u}\RT)^2     \RT) =     \int_{-\infty}^\infty e^{-\frac{ g_I^{''}(t_0) x^2}{4}} \,dx.
		\EQNY
		Hence in view of \ccj{\nelem{lem:HTT}}, letting $L\to\IF$ we obtain
		\begin{eqnarray*}
		\cc{\lim_{u\to\IF}	R_T(u)} \le\frac{1}{t_0^{m/2}}\H_I(T)
			\int_{-\infty}^\infty e^{-\frac{g_I^{''}(t_0) y^2}{4}}\psi(y)\,dy.\nonumber
		\end{eqnarray*}
		Similarly, we obtain the following lower bound
		\BQNY
		\cc{R_T(u)}\ge  \frac{1}{(t_0+\vn)^{m/2}}   F_3(L,Q,\vn,T, u),
		\EQNY
		where
		\BQNY
		F_3(L,Q,\vn,T, u)&=&\frac{T}{\sqrt{u}}\sum_{-N_u-1\le j\le N_u}
		\exp\LT(-\frac{g_I^{''}(t_0)+\vn }{4}  \LT(\frac{jT}{\sqrt u}\RT)^2     \RT)\\
		&&\times
		\int_{\abs{\abs{\vk x_I}}\le L}
		e^{  \frac{  \x _{I }^\top (\SI_{I I }^{-1}\b_I - \vk\vn_I^{\x_I})}{t_0-\vn(\x_I)}  -\vn} P_{j;u}^-(T,\vk x_I,\vn,Q)\,d\vk{x}_I.
		\EQNY
		Letting $u\to\IF, Q\to\IF, \vn\to 0,L\to\IF$ (in \EH{this} order) and in view of \nelem{lem:HTT}  we obtain
		\begin{eqnarray*}
		\cc{\lim_{u\to\IF}	R_T(u)}  \ge\frac{1}{t_0^{m/2}}\H_I(T)
			\int_{-\infty}^\infty e^{-\frac{g_I^{''}(t_0) y^2}{4}}\psi(y)\,dy \EHb{>0}.\nonumber
		\end{eqnarray*}
		Consequently, the claim follows and the proof is complete. \QED

\e1{{\bf Acknowledgement}: Thanks to Swiss National Science Foundation grant No.200021-166274.
	TR \& KD acknowledge partial support by NCN Grant No 2015/17/B/ST1/01102 (2016-2019).
}

\bibliographystyle{ieeetr} 

 \bibliography{vectProcEKEEKK}
\end{document}